\documentclass[12pt]{amsart}
\usepackage{amsmath,amssymb}
\usepackage{mathptmx}
\usepackage{enumerate}
\usepackage{url}
\usepackage[colorlinks=true,citecolor=blue]{hyperref}

\newtheorem{thm}{Theorem}[section]
\newtheorem{lem}[thm]{Lemma}
\newtheorem{prop}[thm]{Proposition}
\newtheorem{cor}[thm]{Corollary}
\newtheorem{knownthm}[thm]{Theorem}
\newtheorem{ques}{Question}
\theoremstyle{definition}
\newtheorem{exam}[thm]{Example}
\newtheorem{rem}[thm]{Remark}
\numberwithin{equation}{section}

\def\N{\mathbb N}
\def\C{\mathbb C}
\def\Z{\mathbb Z}
\def\U{\mathcal U}
\def\B{\mathcal B}
\newcommand{\HH}{\mathcal H}
\newcommand{\Folner}{F\o lner}
\DeclareMathOperator{\ran}{ran}
\DeclareMathOperator{\spec}{spec}

\DeclareMathOperator{\Rea}{Re}
\hypersetup{pdftitle={Antidiagonal and square averages for amenable groups}}

\begin{document}
\title[Antidiagonal and square averages]{Antidiagonal and square averages for amenable groups}
\date{}

\author{Xiaofang Luo}
\address{School of Mathematics, South China University of Technology,
Guangzhou, 510641, P.R. China}
\email{luoxf29@scut.edu.cn}


\begin{abstract}
We prove that, for commuting unitary representations $U,V$ of a countable amenable group $G$, the averages
$$
D_F=\frac1{|F|}\sum_{g\in F}U_g^*V_g
$$
have a common selfadjoint weak operator limit along all left and right \Folner{} sequences, with strong* convergence along two-sided sequences. This answers Bergelson's question on commuting actions. For any unitary representation $\pi$, we also prove that $A_F=|F|^{-1}\sum_{g\in F}\pi(g^2)$ converges strongly along right \Folner{} sequences to a selfadjoint central operator supported on the finite type~I part of $\pi(G)''$; the adjoints $A_F^*$ converge strongly along left \Folner{} sequences. This settles the convergence question for square averages posed by Bergelson and Rodr\'iguez-Mart\'in. For measure-preserving actions of finitely generated amenable or countable virtually FC groups, we obtain
$$
\lim_n\frac1{|F_n|}\sum_{g\in F_n}\mu(A\cap T_{g^2}A)>0
\qquad\text{whenever }\mu(A)>0.
$$
For each $\varepsilon>0$, the set $\{g:\mu(A\cap T_{g^2}A)>\mu(A)^2-\varepsilon\}$ meets every Bohr neighborhood of the identity in a set of positive lower density, uniformly over \Folner{} sequences. These intersections are also left and right IP$^*$ for virtually FC groups. We construct a unitary representation $\pi$ of a locally finite nilpotent group and a weakly mixing vector $f$ for which $\{g:\Rea\langle\pi(g^2)f,f\rangle>-\varepsilon\}$ has upper Banach density zero for small $\varepsilon>0$. We also give a weakly mixing action of a two-generated solvable group for which the IP$^*$ conclusion fails. The recurrence results yield product sets $BB$ inside difference sets, with $B$ of positive upper density.
\end{abstract}

\maketitle

\section{Introduction}

Multiple recurrence connects the return times of dynamical systems with configurations in sets of positive density. Furstenberg's ergodic proof of Szemer\'edi's theorem~\cite{Furstenberg1977} and the extension to commuting transformations by Furstenberg and Katznelson~\cite{FurstenbergKatznelson1978} established this connection. We study two types of correlations for amenable group actions: those between separate translates under two commuting actions, and those between symmetric translates under one action. We ask whether their \Folner{} averages converge and how large the corresponding return sets are when the acting group is nonabelian.

For commuting invertible measure-preserving transformations $R_1,\ldots,R_k$ of a probability space $(X,\B,\mu)$, the averages
$$
\frac1N\sum_{n=0}^{N-1}\prod_{j=1}^k f_j(R_j^n x),
\qquad f_j\in L^\infty(X,\mu),
$$
converge in $L^2(X,\mu)$. Conze and Lesigne~\cite{ConzeLesigne1984} proved this for $k=2$. For powers of a single ergodic transformation, Host and Kra~\cite{HostKra2005} and Ziegler~\cite{Ziegler2007} identified characteristic factors that are inverse limits of nilsystems: replacing the functions by their conditional expectations on these factors leaves the limit unchanged. Frantzikinakis and Kra~\cite{FrantzikinakisKra2005} proved convergence under ergodicity assumptions on the transformations and their pairwise quotients; Tao~\cite{Tao2008} established the general commuting case.

Bergelson and Leibman~\cite{BergelsonLeibman1996} extended multiple recurrence to polynomial iterates, and Leibman~\cite{Leibman1998} to nilpotent groups of transformations. Walsh~\cite{Walsh2012} proved norm convergence for nilpotent polynomial averages; Zorin-Kranich~\cite{ZorinKranich2016} extended this convergence theory to amenable parameter groups. Here the nilpotence assumption is on the group generated by the transformations. Bergelson and Leibman~\cite{BergelsonLeibman2004} exhibited failures of recurrence and convergence for configurations generated by two elements of solvable groups of exponential growth.

For noncommuting transformations, convergence can depend on the choice of iterates. Frantzikinakis and Host~\cite{FrantzikinakisHost2023} proved $L^2$-convergence for the iterates $n,p(n)$, with $\deg p\geq2$, when the transformation in the linear term has zero entropy. Using polynomial Furstenberg joinings, Huang, Shao and Ye~\cite{HuangShaoYe2026Joinings} extended this theorem to several linear iterates of an ergodic zero-entropy transformation and several nonlinear integer-valued polynomial iterates of another ergodic transformation.

Huang, Shao and Ye~\cite{HuangShaoYe2024Counterexample} constructed pairs of ergodic zero-entropy transformations for which convergence fails for prescribed pairs of polynomials of degrees at least five. Kosloff and Sanadhya~\cite{KosloffSanadhya2026} obtained such examples for every prescribed pair of polynomials of degrees at least two, using a two-dimensional local central limit theorem. Huang, Shao and Ye~\cite{HuangShaoYe2025Counterexample} and, independently, Austin~\cite{Austin2025} also obtained zero-entropy counterexamples when both iterates are $n$.

For minimal topological systems, recurrence is expressed through returns of points to open sets. Huang, Shao and Ye~\cite{HuangShaoYe2024Recurrence} proved simultaneous recurrence on a residual set of points for two minimal homeomorphisms, combining a linear iterate of one with iterates of the other given by nonlinear integer-valued polynomials that vanish at zero. For a single minimal transformation, Glasner, Huang, Shao, Weiss and Ye~\cite{GlasnerHuangShaoWeissYe2025} proved that, after almost one-to-one modifications of the map to the maximal pro-nilfactor, the arithmetic-progression orbit closures through a residual set of diagonal points are unions of fibers of the corresponding product factor map. Our two actions commute with one another; the acting group itself may be nonabelian.

For a countable amenable group $G$, intervals are replaced by \Folner{} sets. Bergelson, McCutcheon and Zhang~\cite{BMZ} established norm convergence and a Roth theorem for the triangular configuration
$$
\frac1{|F_n|}\sum_{g\in F_n}f(T_gx)h(S_gT_gx)
$$
of commuting measure-preserving actions $T,S$; joint ergodicity for this configuration was studied by Bergelson and Rosenblatt~\cite{BergelsonRosenblatt1988}. Zorin-Kranich~\cite{ZorinKranich2016} and Austin~\cite{Austin} proved convergence for successive products of several commuting actions, and Austin also proved the corresponding multiple recurrence theorem. Separate translates lead to the following question.

\begin{ques}[Bergelson {\cite[Question~13, p.~56]{Bergelson}}]\label{ques:bergelson}
For commuting measure-preserving actions $T,S$ of a countable amenable group $G$ on a probability space $(X,\B,\mu)$, does
$$
\lim_{n\to\infty}\frac1{|F_n|}\sum_{g\in F_n}\mu(T_gA\cap S_gA)
$$
exist for every $A\in\B$ and every left \Folner{} sequence $(F_n)$?
\end{ques}

When $G$ is abelian, $g\mapsto T_g^{-1}S_g$ is an action, so the mean ergodic theorem applies. For nonabelian $G$ this map need not be an action, and the change of variables does not reduce the question to triangular averages. A related difficulty appears in
$$
\mu(T_g^{-1}A\cap T_gA)=\mu(A\cap T_{g^2}A).
$$
Here $g^2$ is the square of a group element; on the additive group $\Z$ it is $2g$. The squaring map is a homomorphism for abelian groups, but not in general. For nilpotent groups, convergence follows from polynomial ergodic theory, and finite generation gives a uniform positivity bound.

\begin{knownthm}[{\cite[Theorems~1.27 and~1.29]{BRM2026}}]\label{thm:known-nilpotent}
For a measure-preserving action $T$ of a countable nilpotent group $G$ on a probability space $(X,\B,\mu)$ and $A\in\B$, the limit
$$
\lim_{n\to\infty}\frac1{|F_n|}\sum_{g\in F_n}\mu(A\cap T_{g^2}A)
$$
exists and is independent of the left or right \Folner{} sequence. If $G$ is finitely generated, there is a constant $\lambda_G>0$ for which this limit is at least $\lambda_G\mu(A)^2$, uniformly over all such actions and measurable sets.
\end{knownthm}

Bergelson and Rodr\'iguez-Mart\'in formulated the general convergence problem in unitary representations.

\begin{ques}[Bergelson--Rodr\'iguez-Mart\'in {\cite[Question~8.2]{BRM2026}}]\label{ques:square}
If $\pi$ is a unitary representation of a countable amenable group on a Hilbert space $\HH$, does
$$
\lim_{n\to\infty}\frac1{|F_n|}\sum_{g\in F_n}
\langle\pi(g^2)\xi,\xi\rangle
$$
exist for every $\xi\in\HH$ and every left \Folner{} sequence $(F_n)$?
\end{ques}

We answer both questions affirmatively and establish stronger forms of operator convergence. Hilbert spaces below are complex and need not be separable; $\U(\HH)$ denotes the unitary group. Strong* convergence means strong convergence of the operators and their adjoints. The term \emph{antidiagonal} refers to the pair of parameters $(g^{-1},g)$.

\begin{thm}\label{thm:main}
Suppose that $U,V:G\to\U(\HH)$ are unitary representations of a countable amenable group and
$U_gV_h=V_hU_g$ for all $g,h\in G$. The averages
$$
D_{F_n}=\frac1{|F_n|}\sum_{g\in F_n}U_g^*V_g
$$
converge in the weak operator topology along every left or right \Folner{} sequence, with a common selfadjoint contraction $D_{U,V}$ as limit. Moreover,
$$
V_hD_{U,V}=D_{U,V}U_h\qquad(h\in G).
$$
For every two-sided \Folner{} sequence, $D_{F_n}\to D_{U,V}$ strongly*.
\end{thm}

For Koopman representations, this gives convergence of the averages of $\mu(T_gA\cap S_gB)$ for arbitrary measurable $A,B$ (Corollary~\ref{cor:actions}), answering Question~\ref{ques:bergelson}. The limit can be smaller than $\mu(A)^2$ even when the two actions are individually ergodic (Example~\ref{ex:antidiagonal-below-square}).

Square averages require a separate argument, since $g\mapsto\pi(g)^*$ need not be a representation. Their limit is central in the weakly closed algebra $\pi(G)''$ generated by the representation. For a von Neumann algebra $M$, write $Z(M)$ for its center. Its \emph{finite type~I part} consists of the central summands that are measurable fields of finite matrices; their dimensions need not be bounded. Section~\ref{sec:prelim} recalls the operator-algebra terminology used here.

\begin{thm}\label{thm:square-general}
For a unitary representation $\pi:G\to\U(\HH)$ of a countable amenable group, put
$$
M=\pi(G)'',\qquad
A_F=\frac1{|F|}\sum_{g\in F}\pi(g^2)
$$
for nonempty finite $F\subseteq G$. There is a selfadjoint contraction $Q\in Z(M)$ such that
$$
\begin{aligned}
A_{F_n}&\longrightarrow Q\quad\text{strongly for every right \Folner{} sequence},\\
A_{F_n}^*&\longrightarrow Q\quad\text{strongly for every left \Folner{} sequence}.
\end{aligned}
$$
Consequently, $A_{F_n}\to Q$ weakly along every left or right \Folner{} sequence, and strongly* along every two-sided \Folner{} sequence. The operator $Q$ vanishes outside the finite type~I part of $M$.
\end{thm}

Theorem~\ref{thm:square-general} answers Question~\ref{ques:square}. Proposition~\ref{prop:square-limit} identifies $Q$ by the Frobenius--Schur formula, with
$$
\spec(Q)\subseteq\{0\}\cup\{1/d:d\geq1\}\cup\{-1/(2d):d\geq1\}.
$$
If $M$ is finite, the convergence is strong* along every left or right sequence (Corollary~\ref{thm:square-finite}). Without this hypothesis, one example shows the necessity of the averaging directions in both convergence theorems (Example~\ref{ex:strong}); compare the one-sided mean ergodic results of Bj\"orklund and Fish~\cite{BF}.

For a measure-preserving action, the square-correlation mean
$$
\delta_A=\lim_n\frac1{|F_n|}\sum_{g\in F_n}\mu(A\cap T_{g^2}A)
$$
is therefore independent of the left or right \Folner{} sequence. The remaining positivity problem asks whether $\delta_A>0$ whenever $\mu(A)>0$~\cite[Question~1.35]{BRM2026}. Recurrence above $\mu(A)^2-\varepsilon$, especially within a prescribed Bohr neighborhood, asks for more than positivity of the mean.

Lower bounds close to powers of $\mu(A)$ have long been studied in multiple recurrence. For an ergodic transformation, Bergelson, Host and Kra~\cite{BergelsonHostKra2005} proved that the three-term and four-term correlations exceed $\mu(A)^3-\varepsilon$ and $\mu(A)^4-\varepsilon$, respectively, on syndetic sets, using nilsequence decompositions of the correlations. For triangular configurations, Chu and Zorin-Kranich~\cite{ChuZorinKranich2015} obtained a $\mu(A)^4-\varepsilon$ bound on a syndetic set when the action generated by the two commuting actions is ergodic. Jamneshan and Pan~\cite{JamneshanPan2025} proved uniform syndeticity for successive-product configurations over solvable groups of fixed derived length.

We also consider recurrence within prescribed Bohr neighborhoods and IP-sets. Host and Kra~\cite{HostKra2011} studied Nil-Bohr sets, defined by visits of nilrotations to open sets. Huang, Shao and Ye~\cite{HuangShaoYe2016} characterized Nil-Bohr$_0$ sets, given by returns to neighborhoods of the starting point, through generalized polynomials and used them to characterize higher-order almost automorphy. Here we impose a Bohr restriction coming from a compact group rotation. Furstenberg and Katznelson~\cite{FurstenbergKatznelson1985} established multiple recurrence for IP-systems. Bergelson, Furstenberg and McCutcheon~\cite{BFM1996} proved polynomial IP recurrence. Bergelson and McCutcheon~\cite{BergelsonMcCutcheon2007} proved a triangular multiple recurrence theorem for arbitrary countable groups. Parini~\cite{Parini2024} proved left IP$^*$ conclusions for monochromatic triangular configurations over amenable groups, and for their higher-dimensional versions over FC groups.

Our square-recurrence theorem applies to finitely generated amenable groups and to virtually FC groups. A group is \emph{FC} if each conjugacy class is finite, and \emph{virtually FC} if it has an FC subgroup of finite index. Bohr neighborhoods, syndeticity and the two IP$^*$ conventions are recalled in Section~\ref{sec:prelim}.

\begin{thm}\label{thm:bohr}
Suppose that $G$ is countable amenable and is either finitely generated or virtually FC. For a unitary representation $\pi:G\to\U(\HH)$, denote by $P_2$ the orthogonal projection onto
$$
\bigcap_{g\in G}\ker(\pi(g^2)-I).
$$
Given $\xi\in\HH$, $\varepsilon>0$ and a Bohr neighborhood $B$ of the identity, the set
$$
E=\{g\in B:
\Rea\langle\pi(g^2)\xi,\xi\rangle
>\|P_2\xi\|^2-\varepsilon\}
$$
is left and right syndetic. It has a positive lower density bound independent of the averaging sequence: there is $\eta>0$ such that
$$
\liminf_{n\to\infty}\frac{|E\cap F_n|}{|F_n|}\geq\eta
$$
for every left or right \Folner{} sequence $(F_n)$. If $G$ is virtually FC, then $E$ is also left and right IP$^*$.
\end{thm}

For Koopman representations, $\|P_21_A\|_2^2\geq\mu(A)^2$. Thus Theorem~\ref{thm:bohr} gives the asserted return-set bound and, by nonnegativity of correlations, $\delta_A>0$ for every positive-measure $A$. Corollary~\ref{cor:square-recurrence} extends positivity to every nonnegative, nonzero $L^2$ function. This extends the strict-positivity conclusion of Theorem~\ref{thm:known-nilpotent} to finitely generated amenable groups, but does not give a constant $\lambda_G$ uniform over all actions. The lower density bound may depend on $\xi,\varepsilon$ and $B$; Theorem~\ref{thm:bohr} does not assert that $|E\cap F_n|/|F_n|$ has a limit.

For every countable amenable group, Proposition~\ref{prop:square-limit} also gives
$$
\delta_A\geq\frac{\mu(A)(3\mu(A)-1)}2
$$
(Corollary~\ref{cor:large-sets}). Hence $\mu(A)>1/3$ ensures positivity without additional assumptions on $G$. For general countable amenable groups, strict positivity remains open for sets with $0<\mu(A)\leq1/3$.

Theorem~\ref{thm:bohr} does not extend to all countable amenable groups. The next construction answers \cite[Question~8.3(a),(b)]{BRM2026} negatively. The weakly mixing subspace of a representation is the orthogonal complement of the closed span of its finite-dimensional invariant subspaces.

\begin{thm}\label{thm:no-fc-counterexample}
There exist a countable locally finite group $G$, a measure-preserving action $T$ on a probability space $(X,\B,\mu)$ and a set $A\in\B$ with $\mu(A)=1/2$ such that $G$ is nilpotent of class two and exponent four, its center $W=Z(G)$ has infinite index, and
$$
\mu(A\cap T_{g^2}A)=
\begin{cases}
1/2,&g\in W,\\
1/8,&g\notin W.
\end{cases}
$$
For every $0<\varepsilon<1/8$, the return set
$$
\{g\in G:\mu(A\cap T_{g^2}A)>\mu(A)^2-\varepsilon\}
$$
equals $W$. It has density zero along every left or right \Folner{} sequence and is neither left nor right syndetic, and neither left nor right IP$^*$. For the Koopman representation $\pi$ of this action, $f=1_A-1/2$ belongs to the weakly mixing subspace.
\end{thm}

For the centered vector in the theorem,
$$
\langle\pi(g^2)f,f\rangle=\mu(A\cap T_{g^2}A)-\frac14.
$$
Its return set above $-\delta$ is $W$ for $0<\delta<1/8$, although the uncentered mean is $1/8>0$. An earlier example using the quaternion group appears in \cite[Remark~8.4]{BRM2026}. Here the return set is an infinite-index central subgroup, and therefore has density zero along every left or right \Folner{} sequence.

Finite generation also does not ensure IP$^*$ recurrence. In Theorem~\ref{thm:two-generated-square-ip-counterexample}, we construct a weakly mixing action of a two-generated solvable group of derived length at most three. There is a set $A$ of measure $1/2$ whose square return sets above $\mu(A)^2-\varepsilon$, for $0<\varepsilon<1/8$, have density one along every left and right \Folner{} sequence but are neither left nor right IP$^*$. The whole action is weakly mixing, whereas Theorem~\ref{thm:no-fc-counterexample} asserts this only for a particular vector. A related failure of projection limits for polynomial IP-systems appears in \cite[Section~3]{BFM1996}.

There is also a product-set consequence. Bergelson~\cite[Corollary~3.1.2]{Bergelson1985} proved that $D-D$ contains $B+B$ with $B$ of positive upper density whenever $D\subseteq\Z$ has positive upper Banach density; Bergelson and Ruzsa~\cite{BergelsonRuzsa2009} showed that $B$ can be chosen symmetric, containing zero and with positive asymptotic density, and also obtained polynomial variants. Beiglb\"ock, Bergelson and Fish~\cite{BeiglboeckBergelsonFish2010} proved that products of sets of positive left upper Banach density in countable amenable groups are piecewise Bohr. Combining our positivity results with \cite[Theorem~1.24]{BRM2026} gives, for every prescribed left or right \Folner{} sequence $F$,
$$
d_l^*(D)>0\quad\Longrightarrow\quad
\text{there is }B\subseteq G\text{ with }\overline d_F(B)>0
\text{ and }BB\subseteq DD^{-1}
$$
when $G$ is finitely generated amenable or countable virtually FC. For arbitrary countable amenable $G$, the same conclusion holds if $d_l^*(D)>1/3$; see Corollary~\ref{cor:square-products} for the quantitative statement.

Schur orthogonality gives the finite-dimensional antidiagonal limits, as in Choi~\cite[Proposition~3.3]{Choi2023}, and the Frobenius--Schur formula gives the square limits~\cite[Theorem~3.1 and Lemma~3.2]{KoppWG2011}. Lemma~\ref{lem:finite-type-I-reduction} transfers these calculations to the finite type~I part. On the complement, we prove convergence to zero in the directions stated in Theorems~\ref{thm:main} and~\ref{thm:square-general}, using Lemmas~\ref{lem:weighted-multiplication} and~\ref{lem:opposite}. The matrix norm calculation in the latter is related to Tomiyama's transpose-map theorem~\cite{Tomiyama1983}.

For finitely generated amenable groups, we use Breuillard and Pisier's uniform Jordan theorem~\cite[Proposition~1.1]{BreuillardPisier2017} to obtain a finite-index subgroup on which all unitary representations of a prescribed bounded dimension have abelian images (Lemma~\ref{lem:common-abelian-image}). Squaring is a homomorphism on these images. We combine the mean ergodic theorem with a nonnegative weight supported in the given Bohr neighborhood; Lemma~\ref{lem:character-square} controls the remaining terms. For virtually FC groups, we group terms of an arbitrary sequence into finite products that commute with one another to obtain the IP$^*$ assertion. Sections~\ref{sec:proof} and~\ref{sec:square-general} prove convergence; Section~\ref{sec:square-fc} treats square recurrence and product sets, and Sections~\ref{sec:recurrence-counterexamples} and~\ref{sec:examples} contain the counterexamples.

\section{Preliminaries}\label{sec:prelim}

Acting groups are countable and discrete unless stated otherwise; compact groups are Hausdorff. A sequence $F=(F_n)$ of nonempty finite subsets of $G$ is left \Folner{} if
$$
|hF_n\triangle F_n|/|F_n|\longrightarrow0\qquad(h\in G).
$$
Replacing $hF_n$ by $F_nh$ defines a right \Folner{} sequence; a two-sided sequence satisfies both conditions. A countable group is amenable if it admits a left \Folner{} sequence.

For $E\subseteq G$, write
$$
\underline d_F(E)=\liminf_n\frac{|E\cap F_n|}{|F_n|},
\qquad
\overline d_F(E)=\limsup_n\frac{|E\cap F_n|}{|F_n|}.
$$
Their common value, when they agree, is the density along $F$. The left upper Banach density $d_l^*(E)$ is the supremum of $\overline d_F(E)$ over left \Folner{} sequences; $d_r^*(E)$ uses right sequences. A set is left syndetic if $DE=G$ for some finite $D$, and right syndetic if $ED=G$; see \cite{BHM1998} for these notions of size.

A Bohr neighborhood of the identity contains $\rho^{-1}(O)$ for a homomorphism $\rho:G\to K$ into a compact group and an open identity neighborhood $O\subseteq K$. A set is right thick if it contains a right translate of every finite subset of $G$, and piecewise Bohr if it contains the intersection of a translate of a Bohr neighborhood with a right thick set \cite{BeiglboeckBergelsonFish2010}. For a sequence $(g_n)$,
$$
\begin{aligned}
\operatorname{FP}_{\uparrow}(g_n)&=\{g_{i_1}\cdots g_{i_k}:k\geq1,\ i_1<\cdots<i_k\},\\
\operatorname{FP}_{\downarrow}(g_n)&=\{g_{i_k}\cdots g_{i_1}:k\geq1,\ i_1<\cdots<i_k\}.
\end{aligned}
$$
Right IP$^*$ sets meet every $\operatorname{FP}_{\uparrow}(g_n)$; left IP$^*$ sets meet every $\operatorname{FP}_{\downarrow}(g_n)$. No injectivity is required of $(g_n)$.

Actions on probability spaces are left actions, with Koopman convention $\pi(g)f=f\circ T_{g^{-1}}$. Two actions commute if $T_gS_h=S_hT_g$ for all $g,h$. The weakly mixing subspace of a unitary representation is the orthogonal complement of the closed span of its finite-dimensional invariant subspaces. An action is weakly mixing if its diagonal action is ergodic, equivalently if the Koopman representation has no nonzero finite-dimensional invariant subspace in $L^2_0=\{f\in L^2:\int f\,d\mu=0\}$.

We will repeatedly use the following two averaging facts.
\begin{lem}\label{lem:subgroup-density}\label{lem:infinite-index-density}
For a subgroup $L$ of a countable amenable group $G$ and every left
or right \Folner{} sequence $(F_n)$,
$$
\frac{|F_n\cap L|}{|F_n|}\longrightarrow [G:L]^{-1},
\qquad \infty^{-1}:=0.
$$
If $[G:L]<\infty$, then $(F_n\cap L)$ is eventually nonempty and
is a \Folner{} sequence in $L$ on the same side.
\end{lem}

\begin{proof}
For a left \Folner{} sequence and distinct cosets $h_1L,\ldots,h_rL$,
$$
\bigl||F_n\cap h_iL|-|F_n\cap L|\bigr|
\leq |h_i^{-1}F_n\triangle F_n|=o(|F_n|).
$$
Summing gives $\limsup_n|F_n\cap L|/|F_n|\leq1/r$. If the index
is infinite, let $r\to\infty$; if it is finite, sum over all cosets
to obtain the stated limit. Inversion gives the right-sided result.
For finite index, $H_n=F_n\cap L$ is eventually nonempty, and
$$
\frac{|hH_n\triangle H_n|}{|H_n|}
\leq\frac{|F_n|}{|H_n|}
\frac{|hF_n\triangle F_n|}{|F_n|}\longrightarrow0\qquad(h\in L).
$$
The same estimate with $H_nh$ applies on the right.
\end{proof}

\begin{lem}\label{lem:compact-average}
Suppose that $G$ is a discrete amenable group, not necessarily countable, and $\rho:G\to K$ is a homomorphism with dense image in a compact group. For every $f\in C(K)$ and every left or right \Folner{} sequence, or net, $(F_i)$,
$$
\frac1{|F_i|}\sum_{g\in F_i}f(\rho(g))
\longrightarrow\int_K f\,dm_K,
$$
where $m_K$ is normalized Haar measure.
\end{lem}

\begin{proof}
Write $\mu_i=|F_i|^{-1}\sum_{g\in F_i}\delta_{\rho(g)}$. For a left \Folner{} net and $h\in G$,
$$
\left|\int_K f(\rho(h)k)\,d\mu_i(k)-\int_K f(k)\,d\mu_i(k)\right|
\leq\|f\|_\infty\frac{|hF_i\triangle F_i|}{|F_i|}\longrightarrow0.
$$
Every weak* cluster point of $(\mu_i)$ is therefore invariant under $\rho(G)$, hence under $K$ by density and uniform continuity. Uniqueness of Haar probability measure gives $\mu_i\to m_K$. Replacing $f(\rho(h)k)$ and $hF_i$ by $f(k\rho(h))$ and $F_ih$ proves the right-sided assertion.
\end{proof}

Hilbert spaces are complex, with inner products linear in the first variable, and need not be separable. We write $B(\HH)$ for the bounded linear operators and $\U(\HH)$ for the unitary operators. For a bounded sequence $(T_n)$,
$$
\begin{aligned}
T_n\to T\text{ weakly}
&\quad\Longleftrightarrow\quad
\langle T_n\xi,\eta\rangle\to\langle T\xi,\eta\rangle
\quad(\xi,\eta\in\HH),\\
T_n\to T\text{ strongly}
&\quad\Longleftrightarrow\quad
\|(T_n-T)\xi\|\to0\quad(\xi\in\HH).
\end{aligned}
$$
Strong* convergence also requires $T_n^*\to T^*$ strongly. For Koopman representations, weak convergence concerns all $L^2$ correlations, while strong convergence is convergence in $L^2$ for each function.

A von Neumann algebra is a unital $*$-subalgebra of $B(\HH)$ closed in the weak operator topology. For a set $\mathcal S$ of operators,
$$
\mathcal S'=\{a\in B(\HH):as=sa\text{ for every }s\in\mathcal S\}.
$$
The double commutant $M=\pi(G)''$ is the von Neumann algebra generated by the operators $\pi(g)$. Its center is $Z(M)=M\cap M'$. A central projection satisfies $p=p^*=p^2\in Z(M)$ and gives the common reducing decomposition
$$
\HH=p\HH\oplus(1-p)\HH,\qquad
x=x|_{p\HH}\oplus x|_{(1-p)\HH}\quad(x\in M).
$$
For selfadjoint $x\in M$, $1_E(x)$ denotes its spectral projection on a Borel set $E\subseteq\mathbb R$; in particular, $1_{\{a\}}(x)$ projects onto $\ker(x-aI)$.

The algebra $M$ is finite if $v^*v=I$ for $v\in M$ implies $vv^*=I$. A homogeneous type~$I_d$ algebra has the form $M_d(Z)$ with $Z$ abelian. On a separable Hilbert space it acts as a measurable field of $d\times d$ matrices, with possible multiplicities. This matrix description does not mean that the representation is a sum of finite-dimensional invariant subspaces; weakly mixing vectors may occur in it.

Matrix units $(e_{ij})_{i,j=1}^r$ with unit $p$ satisfy
$$
e_{ij}e_{k\ell}=\delta_{jk}e_{i\ell},\qquad
e_{ij}^*=e_{ji},\qquad \sum_i e_{ii}=p.
$$
The projection $p$ need not be the identity of $M$.

\begin{lem}\label{lem:types}
There are pairwise orthogonal central projections $z_0,z_1,\ldots$ with strong sum $1$ such that $z_dM$ is homogeneous of type~$I_d$ for $d\geq1$, while $z_0M$ has no nonzero finite type~I central summand. If $z_0\ne0$, then $z_0M$ contains a copy of $M_{2^k}(\C)$ with unit $z_0$ for every $k\geq1$.
\end{lem}

\begin{proof}
The central type decomposition gives the projections $z_d$ \cite[Theorems~5.3.2, 5.3.17 and~5.3.19]{Nelson}. On $z_0M$, repeated use of the halving lemma \cite[Sections~5.2--5.3]{Nelson} gives, for $r=2^k$, orthogonal projections $q_i$ and operators $v_i$ such that
$$
\sum_{i=1}^r q_i=z_0,\qquad
v_i^*v_i=q_1,\qquad v_iv_i^*=q_i.
$$
Then $e_{ij}=v_iv_j^*$ satisfy the matrix unit relations and have unit $z_0$.
\end{proof}

We use $z_d(M)$ when the algebra needs to be specified, and write
$z_{\mathrm I}(M)=\sum_{d\geq1}z_d(M)$ for its finite type~I part. All sums of central projections are strong operator sums; the finite dimensions $d$ need not be bounded.

A normal linear map between von Neumann algebras is an ultraweakly continuous linear map. On bounded sets, ultraweak convergence agrees with weak operator convergence. In particular, normal maps preserve bounded weak operator limits. For a normal $*$-homomorphism $\varphi:N\to M$, there is a central projection $q\in N$ such that
$$
\ker\varphi=(1-q)N,\qquad qN\cong\varphi(N).
$$
Thus restrictions to reducing subspaces and other normal quotients preserve finite type~I algebras. We use these standard facts in the form given in \cite{Nelson,Takesaki}.

\begin{lem}\label{lem:finite-type-I-reduction}
Suppose that $\Gamma$ is countable and $\pi:\Gamma\to\U(\HH)$ generates a finite type~I algebra. For finitely supported functions $c_n:\Gamma\to\C$, put
$$
T_n=\sum_{g\in\Gamma}c_n(g)\pi(g),\qquad
C=\sup_n\sum_{g\in\Gamma}|c_n(g)|<\infty.
$$
If $\sum_g c_n(g)\rho(g)$ converges in norm for every finite-dimensional unitary representation $\rho$ of $\Gamma$, then $T_n$ converges strongly* to an element of $\pi(\Gamma)''$. Two such sequences have the same limit whenever their limits agree in every finite-dimensional unitary representation.
\end{lem}

\begin{proof}
For separable $\HH$, the central decomposition~\cite{Takesaki} of $M=\pi(\Gamma)''$ gives
$$
\HH=\int_Y^\oplus(\C^{d(y)}\otimes K_y)\,d\nu(y),\qquad
\pi(g)=\int_Y^\oplus(\pi_y(g)\otimes I_{K_y})\,d\nu(y),
$$
where $d(y)<\infty$ and $K_y$ records the multiplicity. Countability of $\Gamma$ puts all relations
$$
\pi_y(e)=I,\qquad \pi_y(g)^*=\pi_y(g^{-1}),\qquad
\pi_y(gh)=\pi_y(g)\pi_y(h)
$$
on one common conull set. There each $\pi_y$ is a unitary representation and generates $M_{d(y)}(\C)$.

By hypothesis, the matrices $T_n(y)=\sum_g c_n(g)\pi_y(g)$
converge in norm to a measurable matrix field $T(y)$, with $\|T_n(y)\|,\|T(y)\|\leq C$. The corresponding operator $T\in M$ satisfies
$$
\|(T_n-T)\xi\|^2
=\int_Y\|((T_n(y)-T(y))\otimes I_{K_y})\xi_y\|^2\,d\nu(y)
\longrightarrow0.
$$
The integrand is bounded by $4C^2\|\xi_y\|^2$, so dominated convergence applies. The adjoint matrix fields converge in the same norm, so the same calculation applies to $T_n^*$.

For arbitrary $\HH$, the cyclic subspace
$\HH_\xi=\overline{\operatorname{span}}\{\pi(g)\xi:g\in\Gamma\}$
is separable and reducing, and the restricted algebra is a normal quotient of $M$. The preceding argument therefore gives limits of $T_n\xi$ and $T_n^*\xi$ for every $\xi$. The operators $T_n$ are uniformly bounded, so $T\xi=\lim_nT_n\xi$ defines a bounded linear operator. Moreover,
$$
\langle T\xi,\eta\rangle
=\lim_n\langle T_n\xi,\eta\rangle
=\left\langle\xi,\lim_n T_n^*\eta\right\rangle.
$$
Hence $T_n^*\eta\to T^*\eta$, and $T\in M$ by strong closedness. Applying the result to the difference of two coefficient sequences proves the final assertion.
\end{proof}

We also use standard bidual extension facts~\cite{Takesaki}. A concrete $C^*$-algebra is a norm-closed selfadjoint algebra of operators. Its representations below are unital and preserve adjoints. For a unital $C^*$-algebra $A$, its bidual $A^{**}$ is a von Neumann algebra. Every bounded linear map $L:A\to M$ has a normal extension $\widetilde L:A^{**}\to M$ of the same norm. If every matrix map $[a_{ij}]\mapsto[L(a_{ij})]$ is contractive, $L$ is called completely contractive; so is $\widetilde L$. A unital $*$-representation extends to a normal $*$-homomorphism onto $L(A)''$.

\section{Antidiagonal averages}\label{sec:proof}

For finite-dimensional $U,V$, write
$K=\overline{\{U_g\oplus V_g:g\in G\}}$. Its coordinate
representations $\alpha,\beta$ commute, so
Lemma~\ref{lem:compact-average} gives
\begin{equation}\label{eq:op-finite-dimensional-limit}
D_{F_n}\longrightarrow
\int_K\alpha(k)^*\beta(k)\,dm_K(k)
\quad\text{in norm}
\end{equation}
for every left or right \Folner{} sequence. The representation
$(k,\ell)\mapsto\alpha(k)\beta(\ell)$ of $K\times K$ is a direct
sum of tensor products of irreducibles. On such a summand, Schur
orthogonality gives
\begin{equation}\label{eq:op-Schur}
\int_K\sigma(k)^*\otimes\tau(k)\,dm_K(k)
=\begin{cases}
d^{-1}\mathsf F_d,&\sigma\simeq\tau,\quad d=\dim\sigma,\\
0,&\sigma\not\simeq\tau,
\end{cases}
\qquad \mathsf F_d(v\otimes w)=w\otimes v.
\end{equation}
Here equivalent irreducibles are identified unitarily. In
particular, every limit in \eqref{eq:op-finite-dimensional-limit}
is selfadjoint and independent of the sequence.

We next estimate the averages on the remaining central summands.
The group algebra consists of finite sums $\sum_g c_gv_g$, with
$v_gv_h=v_{gh}$ and $v_g^*=v_{g^{-1}}$. Its completion for
$$
\left\|\sum_g c_gv_g\right\|
=\sup_\rho\left\|\sum_g c_g\rho(g)\right\|
$$
is the full group $C^*$-algebra $C^*(G)$; the supremum is over all
unitary representations. Every such representation therefore
extends to $C^*(G)$. For a unital $C^*$-algebra $A$, the maximal
tensor product is the completion of finite sums $\sum_i a_i\otimes c_i$
for the norm
$$
\left\|\sum_i a_i\otimes c_i\right\|_{\max}
=\sup_{\sigma,\tau}\left\|\sum_i\sigma(a_i)\tau(c_i)\right\|,
$$
where $\sigma,\tau$ are unital representations with commuting
ranges. Thus this norm bounds the same expression for every
commuting pair of representations.

\begin{lem}\label{lem:weighted-multiplication}
Suppose that $A$ is a unital $C^*$-algebra,
$\pi:A\to B(\HH)$ is a unital representation, and $M=\pi(A)''$.
If $b\in B(\HH)$ is such that
$$
a\otimes c\longmapsto b\pi(a)\pi(c)
$$
extends to a contraction $A\otimes_{\max}A\to B(\HH)$, then
$$
bz_0=0,\qquad \|bz_d\|\leq d^{-1}\quad(d\geq1),
$$
with the projections from Lemma~\ref{lem:types}.
\end{lem}

\begin{proof}
Write $m$ for the given contraction. We use the same symbols
$m,\pi$ for their normal extensions to
$(A\otimes_{\max}A)^{**}$ and $A^{**}$, respectively.
The inclusions $a\mapsto a\otimes1$ and $c\mapsto1\otimes c$
extend to commuting normal $*$-homomorphisms
$j_1,j_2:A^{**}\to(A\otimes_{\max}A)^{**}$. Their defining
identity extends by normality in each variable:
\begin{equation}\label{eq:op-normal-product}
m(j_1(a)j_2(c))=b\pi(a)\pi(c)
\qquad(a,c\in A^{**}).
\end{equation}

There is a central projection $q$ for which the restriction
$\pi:qA^{**}\to M$ is a normal $*$-isomorphism. Hence matrix units
$(e_{ij})_{i,j=1}^r$ with $\sum_i e_{ii}=z\in Z(M)$ lift to
matrix units $(\widehat e_{ij})$ in $qA^{**}$; write
$\widehat z=\sum_i\widehat e_{ii}$. For
$S=\sum_{i,j}j_1(\widehat e_{ij})j_2(\widehat e_{ji})$, the
matrix unit relations and \eqref{eq:op-normal-product} give
$$
\begin{aligned}
S^*&=S,\\
S^2&=\sum_{i,j}j_1(\widehat e_{ii})j_2(\widehat e_{jj})
=j_1(\widehat z)j_2(\widehat z),\\
m(S)&=b\sum_{i,j}e_{ij}e_{ji}=rbz.
\end{aligned}
$$
The product $j_1(\widehat z)j_2(\widehat z)$ is a projection.
Thus $\|S\|\leq1$ and
\begin{equation}\label{eq:op-weighted-bound}
r\|bz\|=\|m(S)\|\leq\|S\|\leq1.
\end{equation}
For $z=z_d$, take $r=d$. For $z=z_0$,
Lemma~\ref{lem:types} permits $r=2^k$ for every $k$;
\eqref{eq:op-weighted-bound} then gives $bz_0=0$.
\end{proof}

\begin{lem}\label{lem:antidiagonal-one-sided-bound}
Write $z_d^V=z_d(V(G)'')$. Along a right \Folner{} sequence,
$$
\limsup_n\|D_{F_n}\xi\|^2
\leq\sum_{d\geq1}d^{-2}\|z_d^V\xi\|^2
\qquad(\xi\in\HH).
$$
In particular, $D_{F_n}$ tends strongly to zero on
$(1-z_{\mathrm I}(V(G)''))\HH$. For left \Folner{} sequences,
the same assertions hold with $U$ in place of $V$.
\end{lem}

\begin{proof}
Put $D_n=D_{F_n}$ and take a weakly convergent subnet
$D_{n_\lambda}^*D_{n_\lambda}\to B$. Then $0\leq B\leq I$,
and every $z_d^V$ commutes with $D_n$, $B$ and $B^{1/2}$.
The contractions
$$
\Phi_F(X)=\frac1{|F|}\sum_{g\in F}V_gXU_g^*
$$
satisfy
$$
\Phi_F(U_hV_k)=U_hD_FV_k,
\qquad D_FV_h=U_hD_{Fh}.
$$
Therefore
\begin{equation}\label{eq:op-translate-estimate}
\|\Phi_{F_n}(U_hV_k)-D_nV_{hk}\|
=\|D_n-D_{F_nh}\|
\leq\frac{|F_nh\triangle F_n|}{|F_n|}\longrightarrow0.
\end{equation}
For finitely supported coefficients $(c_{h,k})$, summing
\eqref{eq:op-translate-estimate} and using $\|\Phi_{F_n}\|\leq1$
gives
$$
\begin{aligned}
\left\|B^{1/2}\Bigl(\sum_{h,k}c_{h,k}V_{hk}\Bigr)\xi\right\|^2
&=\lim_\lambda
\left\|D_{n_\lambda}\Bigl(\sum_{h,k}c_{h,k}V_{hk}\Bigr)\xi\right\|^2\\
&\leq\left\|\sum_{h,k}c_{h,k}U_hV_k\right\|^2\|\xi\|^2\\
&\leq\left\|\sum_{h,k}c_{h,k}v_h\otimes v_k\right\|_{\max}^2\|\xi\|^2.
\end{aligned}
$$
Thus $v_h\otimes v_k\mapsto B^{1/2}V_{hk}$ extends to a
contraction on $C^*(G)\otimes_{\max}C^*(G)$.
Lemma~\ref{lem:weighted-multiplication} applies to
$v_g\mapsto V_g$ and $b=B^{1/2}$, and yields
$$
B^{1/2}z_0^V=0,\qquad \|B^{1/2}z_d^V\|\leq d^{-1}.
$$
Orthogonality of the central projections now gives
$$
\langle B\xi,\xi\rangle
=\sum_{d\geq1}\|B^{1/2}z_d^V\xi\|^2
\leq\sum_{d\geq1}d^{-2}\|z_d^V\xi\|^2.
$$
To obtain the limsup assertion, first take a subsequence realizing
$\limsup_n\|D_n\xi\|^2$, then a weakly convergent subnet of its
second moments. For a left \Folner{} sequence, apply the right-sided
assertion to $(F_n^{-1})$ and use
$D_{F_n}(U,V)=D_{F_n^{-1}}(V,U)$.
\end{proof}

\begin{proof}[Proof of Theorem~\ref{thm:main}]
Write $M_U=U(G)''$, $M_V=V(G)''$ and $N=(M_U\cup M_V)''$.
The projections
$$
p_U=z_{\mathrm I}(M_U),\qquad p_V=z_{\mathrm I}(M_V),\qquad p=p_Up_V
$$
are central in $N$. On $p_{d,e}=z_d(M_U)z_e(M_V)\ne0$, the two
coordinate algebras have types $I_d$ and $I_e$. They are generated
by their centers and matrix units. The two centers are central in
$N$, while products of the two systems of matrix units form a
system of size $de$. Hence
$$
p=\sum_{d,e\geq1}p_{d,e},\qquad p_{d,e}N\text{ has type }I_{de},
$$
and $pN$ is finite type~I.

We apply Lemma~\ref{lem:finite-type-I-reduction} to the unitary
representation $(g,h)\mapsto(U_gV_h)|_{p\HH}$ of $G\times G$ and
the averages over $\{(g^{-1},g):g\in F_n\}$.
In every finite-dimensional representation of $G\times G$, the
coordinate representations commute, so
\eqref{eq:op-finite-dimensional-limit} gives the required common
norm limit. Therefore $D_{F_n}|_{p\HH}$ converges strongly*
for every left or right sequence, independently of that sequence.
The identity $D_F^*=D_{F^{-1}}$ makes this limit selfadjoint.
Extend it by zero on $(1-p)\HH$ and denote it by $D$.

For a right \Folner{} sequence,
Lemma~\ref{lem:antidiagonal-one-sided-bound} applied to $(U,V)$
and $(V,U)$ gives
$$
D_{F_n}(1-p_V)\xi\longrightarrow0,
\qquad D_{F_n}^*(1-p_U)\eta\longrightarrow0.
$$
Since $1-p=(1-p_V)+p_V(1-p_U)$,
$$
\begin{aligned}
|\langle D_{F_n}(1-p)\xi,\eta\rangle|
&\leq\|D_{F_n}(1-p_V)\xi\|\,\|\eta\|\\
&\quad+\|\xi\|\,\|D_{F_n}^*(1-p_U)\eta\|\longrightarrow0.
\end{aligned}
$$
Thus $D_{F_n}\to D$ weakly on $\HH$. For left sequences,
$D_{F_n}=D_{F_n^{-1}}^*$ and $D=D^*$ give the same weak limit.
For a two-sided sequence, the two directions of
Lemma~\ref{lem:antidiagonal-one-sided-bound}, for both ordered
pairs of representations, give
$$
\begin{aligned}
\|D_{F_n}(1-p)\xi\|
&\leq\|D_{F_n}(1-p_U)\xi\|+\|D_{F_n}(1-p_V)\xi\|\longrightarrow0,\\
\|D_{F_n}^*(1-p)\xi\|
&\leq\|D_{F_n}^*(1-p_U)\xi\|+\|D_{F_n}^*(1-p_V)\xi\|\longrightarrow0.
\end{aligned}
$$
Together with convergence on $p\HH$, these inequalities prove
strong* convergence on $\HH$.

Finally, for left \Folner{} sets,
$$
\|V_hD_{F_n}U_h^*-D_{F_n}\|
=\|D_{hF_n}-D_{F_n}\|
\leq\frac{|hF_n\triangle F_n|}{|F_n|}\longrightarrow0.
$$
Taking weak limits gives $V_hD=DU_h$. Also $\|D\|\leq1$, since
each $D_{F_n}$ is a contraction.
\end{proof}

\begin{prop}\label{prop:limit}
For $D=D_{U,V}$ and $N=(U(G)''\cup V(G)'')''$, we have
$D^2\in Z(N)$. For each $d\geq1$, denote by $P_d$ the orthogonal
projection onto $\ker(D^2-d^{-2}I)$ and set $S_d=dDP_d$.
Then $P_d\in Z(N)$, $S_d^*=S_d$, $S_d^2=P_d$, and $S_1=P_1$.
Moreover,
$$
\begin{gathered}
D=\sum_{d\geq1}\frac1dS_d,\qquad
\left\|D-\sum_{d=1}^R\frac1dS_d\right\|\leq\frac1{R+1}\quad(R\geq1),\\
\spec(D)\subseteq\{0,1\}\cup\{\pm1/d:d\geq2\}.
\end{gathered}
$$
For $p=z_{\mathrm I}(U(G)'')z_{\mathrm I}(V(G)'')$, restrict $pN$
to a separable reducing subspace and take its central decomposition.
Every nonzero matrix value of $D$ is unitarily equivalent to
$d^{-1}\mathsf F_d$. The projections $P_d$ and operators $S_d$
do not depend on a choice of matrix units.
\end{prop}

\begin{proof}
Almost every matrix representation of $G\times G$ in the
central decomposition from Lemma~\ref{lem:finite-type-I-reduction}
generates its full matrix algebra, so its coordinate representations
have the form $\sigma\otimes I$ and $I\otimes\tau$ with
$\sigma,\tau$ irreducible. By \eqref{eq:op-Schur}, its limit is
zero unless $\sigma\simeq\tau$, when it equals $d^{-1}\mathsf F_d$.
The identities $\mathsf F_d^*=\mathsf F_d$, $\mathsf F_d^2=I$ and $\mathsf F_1=I$
give $D^2\in Z(N)$, $S_d^*=S_d$, $S_d^2=P_d$, $S_1=P_1$
and the spectral
inclusion. Orthogonality of the eigenspaces gives
$$
D=0\oplus\bigoplus_{d\geq1}d^{-1}S_d,
\qquad
\left\|D-\sum_{d=1}^Rd^{-1}S_d\right\|
=\sup_{d>R}d^{-1}\|S_d\|\leq\frac1{R+1}.
$$
For arbitrary $\HH$, restrict to the separable reducing subspace
generated by each vector. The same identities then hold on all
of $\HH$.
\end{proof}

If $(U\oplus V)(G)''$ has no finite type~I summand, its two
coordinate algebras have none either: each is a normal quotient.
Thus $p=0$ and $D_{U,V}=0$.

The weak assertion also holds for \Folner{} nets in any discrete
amenable group. If weak convergence or independence failed,
two nets, possibly the same, would admit fixed $\xi,\eta$ and
$\varepsilon>0$ such that every pair of tails contains averaging
sets $E,F$ with
$$
|\langle(D_E-D_F)\xi,\eta\rangle|\geq\varepsilon.
$$
Choose such pairs $(E_n,F_n)$ recursively so that their respective one-sided
\Folner{} errors are below $1/n$ on
$$
K_n=\left(\{1\}\cup\bigcup_{j<n}
(E_j\cup E_j^{-1}\cup F_j\cup F_j^{-1})\right)^n.
$$
Every element of $H=\langle\bigcup_n(E_n\cup F_n)\rangle$ belongs
to $K_n$ eventually. Both sequences are therefore \Folner{} in
$H$, contradicting the common limit proved for countable groups.
A bounded weak operator Cauchy net converges; inversion and
translation give selfadjointness and intertwining as above.

For an FC group, the conjugacy class $C_h$ of $h$ is finite and
$$
\begin{aligned}
Fh\setminus F&\subseteq\bigcup_{c\in C_h}(cF\setminus F),\\
|Fh\triangle F|&=2|Fh\setminus F|
\leq2\sum_{c\in C_h}|cF\setminus F|
=\sum_{c\in C_h}|cF\triangle F|.
\end{aligned}
$$
Thus every left \Folner{} sequence is right \Folner{}; inversion
gives the converse. Theorems~\ref{thm:main} and
\ref{thm:square-general} then give strong* convergence along
all one-sided \Folner{} sequences in an FC group.

\begin{cor}\label{cor:actions}
For commuting measure-preserving actions $T,S$ of a countable
amenable group on a probability space $(X,\B,\mu)$ and all
$A,B\in\B$, the limit
$$
\lim_n\frac1{|F_n|}\sum_{g\in F_n}\mu(T_gA\cap S_gB)
$$
exists and is independent of the left or right \Folner{} sequence.
\end{cor}

\begin{proof}
For $U_gf=f\circ T_{g^{-1}}$ and $V_gf=f\circ S_{g^{-1}}$,
Theorem~\ref{thm:main} gives
$$
\begin{aligned}
\mu(T_gA\cap S_gB)
&=\langle V_g1_B,U_g1_A\rangle
=\langle U_g^*V_g1_B,1_A\rangle,\\
\frac1{|F_n|}\sum_{g\in F_n}\mu(T_gA\cap S_gB)
&=\langle D_{F_n}1_B,1_A\rangle
\longrightarrow\langle D_{U,V}1_B,1_A\rangle.
\end{aligned}
$$
Taking $B=A$ answers Bergelson's question.
\end{proof}

\section{Square averages}\label{sec:square-general}

Write $u_g=\pi(g)$, $M=\pi(G)''$ and
$z_{\mathrm I}=\sum_{d\geq1}z_d=1-z_0$, with the central
projections from Lemma~\ref{lem:types}. We first determine the
limit on $z_{\mathrm I}\HH$, then prove strong convergence to zero
on $z_0\HH$ along right \Folner{} sequences.

\begin{lem}[{\cite[Theorem~3.1 and Lemma~3.2]{KoppWG2011}}]\label{lem:FS}
For an irreducible continuous unitary representation
$\rho:K\to\U(V)$ of a compact group, with $\dim V=d$,
$$
\int_K\rho(k^2)\,dm_K(k)=\frac{\nu_\rho}{d}I_V,
\qquad \nu_\rho\in\{0,1,-1\}.
$$
Here $\nu_\rho=\int_K\operatorname{tr}(\rho(k^2))\,dm_K(k)$
is the Frobenius--Schur indicator. It can be negative only when
$d$ is even.
\end{lem}

\begin{lem}\label{lem:finiteI}
If $M$ is finite type~I, the square averages converge strongly*
along every left or right \Folner{} sequence to the same central
selfadjoint operator.
\end{lem}

\begin{proof}
For a finite-dimensional unitary representation $\rho$ of $G$,
Lemma~\ref{lem:compact-average} applied entrywise on
$K=\overline{\rho(G)}$ gives
\begin{equation}\label{eq:sq-finite-dimensional-mean}
\frac1{|F_n|}\sum_{g\in F_n}\rho(g^2)
\longrightarrow\int_K k^2\,dm_K(k)
\quad\text{in norm}.
\end{equation}
Every square average is a finite linear combination of the
operators $\pi(h)$ with nonnegative coefficients summing to one.
Thus \eqref{eq:sq-finite-dimensional-mean} and
Lemma~\ref{lem:finite-type-I-reduction} give strong* convergence
and independence of the \Folner{} sequence.

In the central decomposition on a separable reducing subspace,
almost every matrix representation generates its full matrix
algebra and is irreducible. Lemma~\ref{lem:FS} makes its limit a
real scalar matrix. Hence the limit is central and selfadjoint.
The separable reducing subspaces in
Lemma~\ref{lem:finite-type-I-reduction} give these identities on
all of $\HH$.
\end{proof}

For a von Neumann algebra $N$, the opposite algebra has product
$x^{\mathrm{op}}y^{\mathrm{op}}=(yx)^{\mathrm{op}}$ and the same
linear space, norm and adjoint as $N$. For $N\subseteq B(\HH)$,
matrices in $M_r(N)$ act on $\HH^r$; their operator norms are
the matrix norms used below.

\begin{lem}\label{lem:opposite}
For $a\in N$, complete contractivity of either map
$$
N^{\mathrm{op}}\longrightarrow N,\qquad
x^{\mathrm{op}}\longmapsto xa
\quad\text{or}\quad x^{\mathrm{op}}\longmapsto ax
$$
forces $a$ to vanish outside the finite type~I part of $N$.
\end{lem}

\begin{proof}
Work on the complementary central summand, first for
$x^{\mathrm{op}}\mapsto xa$. If $a\ne0$, choose $\varepsilon>0$
with $p=1_{[\varepsilon^2,\infty)}(aa^*)\ne0$. After restricting
to a central summand where $p\ne0$, there are matrix units $(e_{ij})_{i,j=1}^r$ of arbitrarily large size with
$0\ne e_{11}\leq p$. On a type~II or III summand, apply
Lemma~\ref{lem:types} to $pNp$ and take $r=2^k$. On an infinite
homogeneous type~I summand, take a nonzero abelian projection
below $p$ and $r$ orthogonal equivalent projections in $N$.
Their sum need not be $1_N$.

Put $p_0=\sum_i e_{ii}$, $X=[e_{ij}^{\mathrm{op}}]$ and
$Y=[e_{ij}a]$. The matrix unit relations give
$$
\begin{gathered}
\relax [e_{ij}]=[e_{i1}]_i[e_{1j}]_j\geq0,\qquad
[e_{ij}]^2=r[e_{ij}]\ne0,\\
\relax [e_{ji}]^*=[e_{ji}],\qquad
[e_{ji}]^2=I_r\otimes p_0\ne0.
\end{gathered}
$$
Matrix transposition identifies $M_r(N^{\mathrm{op}})$ isometrically
with $M_r(N)^{\mathrm{op}}$. Hence
\begin{equation}\label{eq:sq-opposite-norms}
\|X\|=\|[e_{ji}]\|=1,\qquad [e_{ij}]\geq0,\qquad
\|[e_{ij}]\|=r.
\end{equation}
Since $aa^*\geq\varepsilon^2e_{11}$,
$$
YY^*=\left[\sum_k e_{ik}aa^*e_{kj}\right]
\geq\varepsilon^2
\left[\sum_k e_{ik}e_{11}e_{kj}\right]
=\varepsilon^2[e_{ij}].
$$
By \eqref{eq:sq-opposite-norms}, complete contractivity implies
$$
\varepsilon^2r\leq\|Y\|^2\leq\|X\|^2=1,
$$
a contradiction for large $r$. For $x^{\mathrm{op}}\mapsto ax$,
the assumed matrix inequality gives
$$
\|[x_{ij}a^*]\|
=\|[ax_{ji}^*]\|
\leq\|[(x_{ji}^*)^{\mathrm{op}}]\|
=\|[x_{ij}^{\mathrm{op}}]\|.
$$
Thus $x^{\mathrm{op}}\mapsto xa^*$ is completely contractive,
and the first case applies to $a^*$.
\end{proof}

\begin{proof}[Proof of Theorem~\ref{thm:square-general}]
Extend the limit $Q$ from Lemma~\ref{lem:finiteI} by zero on
$z_0\HH$. Fix a right \Folner{} sequence and write $A_n=A_{F_n}$.
We also use $\pi$ for the representation of
$\mathcal A=C^*(G)$ given by $\pi(v_g)=u_g$. The map $\rho:\mathcal A\to M^{\mathrm{op}}$ given by
$\rho(v_g)=(u_g^*)^{\mathrm{op}}$ is a $*$-representation, since
$$
(u_g^*)^{\mathrm{op}}(u_h^*)^{\mathrm{op}}
=(u_h^*u_g^*)^{\mathrm{op}}=(u_{gh}^*)^{\mathrm{op}}.
$$
The maps
$$
T_n(b)=\frac1{|F_n|}\sum_{g\in F_n}u_g\pi(b)u_g
$$
are completely contractive: at each matrix size, each summand
is multiplication on both sides by unitary matrices. Right
translation gives
$$
T_n(v_h)u_h=A_{F_nh},\qquad
\|T_n(v_h)-A_nu_h^*\|
\leq\frac{|F_nh\triangle F_n|}{|F_n|}\longrightarrow0.
$$
Since $\|T_n(b)\|,\|A_nx\|\leq\|b\|$ when
$\rho(b)=x^{\mathrm{op}}$, linearity and density give
\begin{equation}\label{eq:asymptotic-map}
\|T_n(b)-A_nx\|\longrightarrow0
\qquad\bigl(b\in\mathcal A,\ \rho(b)=x^{\mathrm{op}}\bigr).
\end{equation}

Take a weakly convergent subnet $A_{n_\lambda}^*A_{n_\lambda}\to B$;
then $B\in M$ and $0\leq B\leq I$. For $[b_{ij}]\in M_r(\mathcal A)$
with $\rho(b_{ij})=x_{ij}^{\mathrm{op}}$,
\eqref{eq:asymptotic-map} and the matrix norm bound for $T_n$ give,
for every $\xi\in\HH^r$,
$$
\begin{aligned}
\|[B^{1/2}x_{ij}]\xi\|^2
&=\lim_\lambda\|(I_r\otimes A_{n_\lambda})[x_{ij}]\xi\|^2\\
&\leq\|[b_{ij}]\|^2\|\xi\|^2.
\end{aligned}
$$
Thus $L:\mathcal A\to M$, defined by $L(b)=B^{1/2}x$ when
$\rho(b)=x^{\mathrm{op}}$, is completely contractive.

Keep the symbols $L,\rho$ for their normal extensions to
$\mathcal A^{**}$. Since $\rho(\mathcal A)''=M^{\mathrm{op}}$,
a central projection $e\in\mathcal A^{**}$ makes
$\rho:e\mathcal A^{**}\to M^{\mathrm{op}}$
a normal $*$-isomorphism. Both sides of $L(y)=B^{1/2}x$, with
$\rho(y)=x^{\mathrm{op}}$, are normal in $y$, so this identity
holds on $\mathcal A^{**}$. For $[x_{ij}^{\mathrm{op}}]\in M_r(M^{\mathrm{op}})$,
lift its entries to
$y_{ij}\in e\mathcal A^{**}$. The isomorphism preserves matrix
norms, so
$$
\|[B^{1/2}x_{ij}]\|
=\|[L(y_{ij})]\|
\leq\|[y_{ij}]\|
=\|[x_{ij}^{\mathrm{op}}]\|.
$$
Lemma~\ref{lem:opposite} now applies to
$x^{\mathrm{op}}\mapsto B^{1/2}x$ and gives
\begin{equation}\label{eq:sq-cluster-support}
z_0B=0.
\end{equation}
Since \eqref{eq:sq-cluster-support} holds for every weak cluster
point, $z_0A_n^*A_n\to0$ weakly. Hence
$$
\|A_nz_0\xi\|^2
=\langle A_n^*A_nz_0\xi,z_0\xi\rangle\longrightarrow0.
$$
Together with Lemma~\ref{lem:finiteI}, this proves $A_n\to Q$
strongly on $\HH$.

For a left \Folner{} sequence, inversion gives
$A_{F_n}^*=A_{F_n^{-1}}\to Q$ strongly. Since $Q=Q^*\in Z(M)$,
these two assertions give weak convergence for both sides and
strong* convergence for two-sided sequences.
\end{proof}

\begin{prop}\label{prop:square-limit}
After restricting $z_dM$ to a separable reducing subspace of
$z_d\HH$ and taking the central decomposition, the limit in
Theorem~\ref{thm:square-general} is
$$
Q(\omega)=\frac{\nu_\omega}{d}I_d,\qquad
\nu_\omega\in\{0,1,-1\},
$$
where $\nu_\omega$ is the Frobenius--Schur indicator of the
irreducible representation of the compact closure of $\pi(G)$
on almost every matrix factor. A negative indicator requires
even $d$. Moreover,
$$
\left\|Q-\sum_{d=1}^Rz_dQ\right\|\leq\frac1{R+1},\qquad
\spec(Q)\subseteq\{0\}\cup\{1/d:d\geq1\}
\cup\{-1/(2d):d\geq1\}.
$$
For the projection $P_2$ onto the vectors fixed by every
$\pi(g^2)$,
\begin{equation}\label{eq:square-lower}
P_2=1_{\{1\}}(Q),\qquad
Q\geq\frac32P_2-\frac12I.
\end{equation}
\end{prop}

\begin{proof}
Equation~\eqref{eq:sq-finite-dimensional-mean} and
Lemma~\ref{lem:FS} give the stated matrix formula. Orthogonality
of the central projections yields
$$
\|z_dQ\|\leq d^{-1},\qquad
\left\|Q-\sum_{d=1}^Rz_dQ\right\|
=\sup_{d>R}\|z_dQ\|\leq(R+1)^{-1}.
$$
The same matrix formula gives the spectral inclusion on each
separable reducing subspace, and hence on $\HH$.

The subgroup $G^{(2)}=\langle g^2:g\in G\rangle$ is normal and
amenable. For a right \Folner{} sequence $(E_n)$ in this subgroup,
the mean ergodic theorem gives
$$
\frac1{|E_n|}\sum_{h\in E_n}\pi(h)\longrightarrow P_2
\quad\text{strongly},
$$
so $P_2\in M$. Normality of $G^{(2)}$ implies
$\pi(g)P_2\HH=P_2\HH$ for every $g$, and thus $P_2\in Z(M)$.
The quotient $G/G^{(2)}$ has exponent two and is abelian.
Consequently $P_2M$ is abelian and $P_2\leq z_1$.
On $z_1\HH$, the map $g\mapsto z_1u_g^2$ is a unitary
representation with fixed space $P_2\HH$. A second application
of the mean ergodic theorem gives
$$
z_1Q=P_2,\qquad Q=P_2+(1-z_1)Q,\qquad
\|(1-z_1)Q\|\leq\frac12.
$$
The last bound excludes the eigenvalue $1$ on $(1-z_1)\HH$, so
$1_{\{1\}}(Q)=P_2$. Moreover,
$$
Q\geq P_2-\frac12(1-z_1)
\geq P_2-\frac12(1-P_2)
=\frac32P_2-\frac12I.
$$
\end{proof}

\begin{cor}\label{cor:large-sets}
For a measure-preserving action of a countable amenable group,
the common square-correlation limit of a measurable set $A$
with $\mu(A)=\alpha$ satisfies
$$
\delta_A\geq
\frac32\|\mathbb E(1_A\mid\mathcal I_{G^{(2)}})\|_2^2
-\frac12\alpha
\geq\frac{\alpha(3\alpha-1)}2,
$$
where $\mathcal I_{G^{(2)}}$ is the sigma-algebra of sets
invariant modulo null sets under $G^{(2)}$.
In particular, $\delta_A>0$ if $\alpha>1/3$.
\end{cor}

\begin{proof}
For the Koopman representation,
$$
P_21_A=\mathbb E(1_A\mid\mathcal I_{G^{(2)}}),\qquad
\|P_21_A\|_2^2\geq\left(\int P_21_A\,d\mu\right)^2=\alpha^2.
$$
Theorem~\ref{thm:square-general} and \eqref{eq:square-lower} give
$$
\delta_A=\langle Q1_A,1_A\rangle
\geq\frac32\|P_21_A\|_2^2-\frac12\|1_A\|_2^2
\geq\frac{\alpha(3\alpha-1)}2.
$$
\end{proof}

The bound is attained at $\alpha=1/2$: for left translation of
$Q_8$ on itself and $A=\{1,i,j,k\}$, the two elements with square
$1$ contribute $1/2$ and the other six contribute zero, giving
$\delta_A=1/8$; see \cite[Remark~1.26, equation~(14)]{BRM2026}.
This does not establish optimality of the threshold $1/3$.

\begin{cor}\label{thm:square-finite}
Under the hypotheses of Theorem~\ref{thm:square-general}, if
$\pi(G)''$ is finite, the square averages converge strongly*
to $Q$ along every left or right \Folner{} sequence.
\end{cor}

\begin{proof}
A finite von Neumann algebra has a faithful normal center-valued
trace $\operatorname{Tr}_Z:M\to Z(M)$
\cite[Theorem~6.1.12]{Nelson}. The properties needed here are
$$
\operatorname{Tr}_Z(xy)=\operatorname{Tr}_Z(yx),\qquad
x\geq0,\ \operatorname{Tr}_Z(x)=0\ \Longrightarrow\ x=0.
$$
For a right \Folner{} sequence, Theorem~\ref{thm:square-general}
gives $b_n=A_{F_n}-Q\to0$ strongly, with $\|b_n\|\leq2$. Thus
$$
\|b_n^*b_n\xi\|\leq2\|b_n\xi\|\longrightarrow0,
\qquad
\operatorname{Tr}_Z(b_nb_n^*)
=\operatorname{Tr}_Z(b_n^*b_n)\longrightarrow0
\quad\text{ultraweakly}.
$$
Every ultraweak cluster point $C$ of $b_nb_n^*$ is positive and,
by normality, satisfies $\operatorname{Tr}_Z(C)=0$. Faithfulness
gives $C=0$. Therefore
$$
b_nb_n^*\longrightarrow0\quad\text{ultraweakly},\qquad
\|b_n^*\xi\|^2=\langle b_nb_n^*\xi,\xi\rangle\longrightarrow0.
$$
This proves strong* convergence on the right. Inversion gives
the left \Folner{} assertion.
\end{proof}

\section{Square recurrence}\label{sec:square-fc}

For amenable groups, finite generation provides a common finite-index subgroup on which
all unitary representations of bounded dimension have abelian images.
The required uniform bound comes from
\cite[Proposition~1.1]{BreuillardPisier2017}: every amenable discrete
subgroup of $U(d)$ has an abelian normal subgroup of index at most
$J(d)<\infty$, without closedness or finite generation assumptions.

\begin{lem}\label{lem:common-abelian-image}
For each integer $R\geq1$, a finitely generated amenable group $G$ has a
finite-index normal subgroup $C_R$ such that $\rho(C_R)$ is abelian
for every unitary representation $\rho:G\to U(d)$ with $d\leq R$.
\end{lem}

\begin{proof}
For an $s$-generated group,
$|\operatorname{Hom}(G,S_m)|\leq(m!)^s$.
Each subgroup of index $m$ occurs as a point stabilizer in such
an action, so there are at most $m(m!)^s$. Hence
$$
j_R=\max_{1\leq d\leq R}J(d),\qquad
C_R=\bigcap_{\substack{L\lhd G\\{}[G:L]\leq j_R}}L
$$
is normal of finite index. For $\rho:G\to U(d)$, $d\leq R$,
\cite[Proposition~1.1]{BreuillardPisier2017} supplies an abelian
normal $A_\rho\lhd\rho(G)$ with
$$
[G:\rho^{-1}(A_\rho)]=[\rho(G):A_\rho]\leq J(d)\leq j_R,
\qquad \rho(C_R)\subseteq A_\rho.
$$
\end{proof}

\begin{lem}\label{lem:fg-finite-index-type-I}
If $C\lhd G$ has finite index and $u:G\to\U(\HH)$ is unitary,
then $N=u(C)''$ and $M=u(G)''$ satisfy
$z_{\mathrm I}(N)=z_{\mathrm I}(M)$ as projections on $\HH$.
\end{lem}

\begin{proof}
We first check that von Neumann subalgebras of finite type~I
algebras remain finite type~I. Suppose $A\subseteq M_d(Z)$, with $Z$
abelian. If $1-\sum_{k=1}^d z_k(A)\ne0$, Lemma~\ref{lem:types}
gives nonzero matrix units $(e_{ij})_{i,j=1}^{d+1}$ in $A$.
Representing $Z$ as $C(\Omega)$, choose $\omega$ with
$e_{11}(\omega)\ne0$. Then
$$
d\geq\operatorname{rank}\sum_{i=1}^{d+1}e_{ii}(\omega)
=(d+1)\operatorname{rank}e_{11}(\omega)\geq d+1,
$$
a contradiction. For $A\subseteq\bigoplus_dM_d(Z_d)$, the normal
coordinate maps $\theta_d$ have central kernels, with
$$
\ker\theta_d=(1-e_d)A\quad(e_d\in Z(A)),\qquad
e_dA\cong\theta_d(A),\qquad
1=\bigvee_de_d\leq z_{\mathrm I}(A).
$$
The equality holds because the coordinate maps separate points.
Central kernels also show that normal quotients preserve finite type~I.

Suppose now that $N$ is finite type~I. For coset representatives
$r_q$, $q\in G/C$, write $k=[G:C]$ and $v_q=u(r_q)$. Define matrices
in $M_k(N)$ by
$$
(\widetilde u_g)_{q,r}=
\begin{cases}v_q^*u_gv_r,&q=gr,\\0,&q\ne gr.\end{cases}
$$
Indeed, $q=gr$ implies $r_q^{-1}gr_r\in C$. Each matrix has one
unitary entry in each row and column. Multiplication gives
$(v_q^*u_gv_r)(v_r^*u_hv_s)=v_q^*u_{gh}v_s$, so
$\widetilde u_g\widetilde u_h=\widetilde u_{gh}$.
The isometry $W:\HH\to\HH^k$ given below has reducing range:
$$
W\xi=(k^{-1/2}v_q^*\xi)_q,\qquad
\widetilde u_gW=Wu_g,\qquad
W^*(\widetilde u(G)'')W=M.
$$
Thus $M$ is a normal quotient of a subalgebra of $M_k(N)$ and is
finite type~I.

In general, put $z=z_{\mathrm I}(N)$ and $p=z_{\mathrm I}(M)$.
Since $C\lhd G$, we have $u_gNu_g^*=N$ and $u_gzu_g^*=z$.
Thus $z\in Z(M)$, and the argument on $z\HH$ gives $z\leq p$. Conversely,
$pN\subseteq pM$ is finite type~I. The normal map $x\mapsto px$
on $N$ has a central kernel, so for some $e\in Z(N)$,
$$
\ker(x\mapsto px)=(1-e)N,\qquad
eN\cong pN,\qquad e\leq z,\qquad p(1-e)=0.
$$
Hence $p(1-z)=p(1-e)(1-z)=0$, which proves $p\leq z$.
\end{proof}

Continuous functions on compact abelian groups will impose the Bohr
restriction. We first treat unitary characters, that is,
homomorphisms into $\mathbb T$.

\begin{lem}\label{lem:character-square}\label{cor:abelian-weight}
Suppose $G$ is countable amenable, $\pi:G\to\U(\HH)$ is unitary,
and $\rho:G\to L$ is a homomorphism into a compact abelian group.
For $w\in C(L)$, the averages
$$
A_F^w=|F|^{-1}\sum_{g\in F}w(\rho(g))\pi(g^2)
$$
converge weakly, along every left or right \Folner{} sequence,
to the same $Q_w\in Z(M)$, where $M=\pi(G)''$, with
$$
\|Q_w\|\leq\|w\|_\infty,\qquad Q_w(1-z_{\mathrm I}(M))=0.
$$
The averages converge strongly to $Q_w$ on the right, their adjoints
converge strongly to $Q_w^*$ on the left, and convergence is strong*
on two-sided sequences.
If $w\circ\rho$ is a unitary character, $Q_w$ is a selfadjoint contraction.
\end{lem}

\begin{proof}
For a unitary character $\chi:G\to\mathbb T$, the group
$\widetilde G=\{(g,z)\in G\times\mathbb T:z^2=\chi(g)\}$ is a
central extension of $G$ by $\{\pm1\}$, hence countable and
amenable. Its representation $\widetilde\pi(g,z)=z\pi(g)$ satisfies
\begin{equation}\label{eq:rec-cover-algebra}
\widetilde\pi(\widetilde G)''=M,\qquad
\widetilde\pi((g,z)^2)=\chi(g)\pi(g^2).
\end{equation}
For $\widetilde F=\{(g,z)\in\widetilde G:g\in F\}$,
$|\widetilde F|=2|F|$. If $\widetilde h=(h,z)\in\widetilde G$, then
$$
\begin{gathered}
|\widetilde h\widetilde F\triangle\widetilde F|=2|hF\triangle F|,
\qquad
|\widetilde F\widetilde h\triangle\widetilde F|=2|Fh\triangle F|,\\
\frac1{|\widetilde F|}\sum_{\widetilde g\in\widetilde F}
\widetilde\pi(\widetilde g^2)
=\frac1{|F|}\sum_{g\in F}\chi(g)\pi(g^2).
\end{gathered}
$$
Thus Theorem~\ref{thm:square-general} and
\eqref{eq:rec-cover-algebra} give all assertions for $\chi$.

By Stone--Weierstrass, choose finite linear combinations $w_j$
of continuous characters of $L$ with
$\|w-w_j\|_\infty\to0$. The first part gives their limits, and
$$
\|A_F^{w_j}-A_F^{w_k}\|\leq\|w_j-w_k\|_\infty,
\qquad
\|Q_{w_j}-Q_{w_k}\|\leq\|w_j-w_k\|_\infty.
$$
Hence $Q_{w_j}\to Q_w$ in norm, with
$Q_w\in Z(M)$, $Q_w(1-z_{\mathrm I}(M))=0$, and
$\|Q_w\|\leq\|w\|_\infty$. Along a right \Folner{} sequence,
$$
\limsup_n\|(A_{F_n}^w-Q_w)\xi\|
\leq2\|w-w_j\|_\infty\|\xi\|\longrightarrow0.
$$
The same norm estimate for adjoints proves the left-sided assertion;
the weak and two-sided strong* assertions follow.
\end{proof}

A function $w:L\to\C$ is positive definite if
$\sum_{i,j}a_i\overline{a_j}w(k_j^{-1}k_i)\geq0$ for all finite
families $a_i\in\C$, $k_i\in L$.

\begin{lem}\label{lem:pd-cutoff}
For a compact abelian group $L$ with normalized Haar measure $m_L$
and an open identity neighborhood $O$,
there is $u\in C(L)$ with $u\geq0$ and $\|u\|_2=1$ such that
$$
w(k)=\langle\lambda(k)u,u\rangle,\qquad
(\lambda(k)u)(x)=u(k^{-1}x)
$$
is continuous and positive definite, and
$$
0\leq w\leq1,\quad w(e)=1,\quad \operatorname{supp}w\subseteq O,\quad
a:=\int_Lw\,dm_L>0.
$$
For every amenable group $C$ and homomorphism $\rho:C\to L$ with dense image,
the mean of $w\circ\rho$ along every left or right \Folner{} sequence
is $a$.
\end{lem}

\begin{proof}
Choose an identity neighborhood $V$ with
$\overline V\,\overline V^{-1}\subseteq O$ and a nonnegative
$u\in C(L)$ supported in $V$, with $\|u\|_2=1$. Then
$$
\begin{aligned}
w(k)&=\int_Lu(k^{-1}x)u(x)\,dm_L(x),\\
\sum_{i,j}a_i\overline{a_j}w(k_j^{-1}k_i)
&=\Bigl\|\sum_i a_i\lambda(k_i)u\Bigr\|_2^2\geq0.
\end{aligned}
$$
The integral formula gives continuity, and Cauchy--Schwarz gives
$$
\begin{gathered}
0\leq w(k)\leq\|\lambda(k)u\|_2\|u\|_2=1,
\qquad w(e)=1,\\
\operatorname{supp}w\subseteq
(\operatorname{supp}u)(\operatorname{supp}u)^{-1}\subseteq O,
\qquad
\int_Lw\,dm_L=\Bigl(\int_Lu\,dm_L\Bigr)^2>0.
\end{gathered}
$$
Apply Lemma~\ref{lem:compact-average} to $w$ for the mean assertion.
\end{proof}

\begin{proof}[Proof of Theorem~\ref{thm:bohr}, finitely generated case]
Since the conclusions pass to supersets, we may shrink $B$.
Finite-dimensional continuous unitary representations separate
points of compact groups, so compactness gives
$\rho:G\to U(m)$ with $B=\rho^{-1}(O)$ for a symmetric open
identity neighborhood $O$.
For $M=\pi(G)''$, choose $R\geq m$ first, with
\begin{equation}\label{eq:bohr-tail}
\begin{gathered}
p_R=\sum_{d=1}^Rz_d(M),\qquad
q_R=z_{\mathrm I}(M)-p_R,\qquad z_0=1-z_{\mathrm I}(M),\\
\|q_R\xi\|^2<\varepsilon/2.
\end{gathered}
\end{equation}
Lemma~\ref{lem:common-abelian-image} gives $C=C_R$ and makes
$L=\overline{\rho(C)}$ abelian. Lemma~\ref{lem:pd-cutoff}, applied
to $O\cap L$, gives $u,w$ with
$$
0\leq w\leq1,\qquad
\{h\in C:w(\rho(h))\ne0\}\subseteq B,
\qquad a=\int_Lw\,dm_L>0.
$$
These choices do not depend on a \Folner{} sequence. For every
right \Folner{} sequence $(H_n)$ in $C$,
\begin{equation}\label{eq:weight-mean}
\frac1{|H_n|}\sum_{h\in H_n}w(\rho(h))\longrightarrow a.
\end{equation}

On a separable reducing subspace of $p_R\HH$, the decomposition
in Lemma~\ref{lem:finite-type-I-reduction} gives matrix
representations $\pi_\omega$ of dimensions at most $R$.
Countability of $G$ lets all group identities hold on one conull set.
Lemma~\ref{lem:common-abelian-image} then gives, for $h,k\in C$,
$$
\pi_\omega(h)\pi_\omega(k)=\pi_\omega(k)\pi_\omega(h),
\qquad
\pi((hk)^2)p_R=\pi(h^2)\pi(k^2)p_R.
$$
For arbitrary $\HH$, apply this to
$\overline{\operatorname{span}}\{\pi(g)\zeta:g\in G\}$ for each
$\zeta\in p_R\HH$. Thus
$h\mapsto\pi(h^2)|_{p_R\HH}$ is a representation of $C$.

Proposition~\ref{prop:square-limit} gives $P_2\leq z_1(M)\leq p_R$.
Set $t=\|P_2\xi\|^2$ and $\eta=(p_R-P_2)\xi$. The identities
$$
\begin{aligned}
\langle\pi(h^2)p_R\xi,p_R\xi\rangle
&=t+\langle\pi(h^2)\eta,\eta\rangle,\\
w(\rho(h))\langle\pi(h^2)\eta,\eta\rangle
&=\langle(\lambda(\rho(h))\otimes\pi(h^2))
(u\otimes\eta),u\otimes\eta\rangle
\end{aligned}
$$
express the weighted second term as a coefficient of a unitary
representation on $L^2(L)\otimes p_R\HH$. Its mean is the squared
norm of the invariant projection of $u\otimes\eta$. The mean
ergodic theorem and \eqref{eq:weight-mean} give
\begin{equation}\label{eq:rec-weighted-low}
\lim_n\frac1{|H_n|}\sum_{h\in H_n}w(\rho(h))
\Rea\langle\pi(h^2)p_R\xi,p_R\xi\rangle\geq at.
\end{equation}
Lemma~\ref{lem:fg-finite-index-type-I} identifies
$z_{\mathrm I}(\pi(C)'')=z_{\mathrm I}(M)$, so
Lemma~\ref{cor:abelian-weight} applied to $\pi|_C$ gives
\begin{equation}\label{eq:weighted-zero}
\frac1{|H_n|}\sum_{h\in H_n}w(\rho(h))
\langle\pi(h^2)z_0\xi,z_0\xi\rangle\longrightarrow0.
\end{equation}
The remaining part satisfies
$\Rea\langle\pi(h^2)q_R\xi,q_R\xi\rangle\geq-\|q_R\xi\|^2$.

Write $c(g)=\Rea\langle\pi(g^2)\xi,\xi\rangle$. Since
$c(g^{-1})=c(g)$, we have $E=E^{-1}$. Orthogonality and
\eqref{eq:rec-weighted-low}, \eqref{eq:weighted-zero},
\eqref{eq:weight-mean} yield
\begin{equation}\label{eq:weighted-surplus}
\liminf_n\frac1{|H_n|}\sum_{h\in H_n}w(\rho(h))(c(h)-t+\varepsilon)
\geq a(\varepsilon-\|q_R\xi\|^2)>a\varepsilon/2,
\end{equation}
where the last inequality is \eqref{eq:bohr-tail}. With
$b=\|\xi\|^2$, we have
$$
w(\rho(h))(c(h)-t+\varepsilon)
\leq(b-t+\varepsilon)1_E(h).
$$
Indeed, outside $E$ either $w(\rho(h))=0$ or
$c(h)\leq t-\varepsilon$, and on $E$ we use $w\leq1$ and $c\leq b$.
Equation~\eqref{eq:weighted-surplus} now gives
\begin{equation}\label{eq:rec-density-in-C}
\liminf_n\frac{|H_n\cap E|}{|H_n|}
\geq\frac{a\varepsilon}{2(b-t+\varepsilon)}.
\end{equation}
For every right \Folner{} sequence $(F_n)$ in $G$,
Lemma~\ref{lem:subgroup-density} makes $H_n=F_n\cap C$ right
\Folner{} in $C$ and $|H_n|/|F_n|\to[G:C]^{-1}$. Thus
\eqref{eq:rec-density-in-C} implies
\begin{equation}\label{eq:uniform-bohr-density}
\liminf_n\frac{|F_n\cap E|}{|F_n|}
\geq\frac{a\varepsilon}{2[G:C](b-t+\varepsilon)}>0.
\end{equation}
Since $E=E^{-1}$, the equality
$|F_n\cap E|=|F_n^{-1}\cap E|$ gives the same left-sided bound.
If $E$ were not right syndetic, choose $x_n\notin EF_n^{-1}$.
Then
$$
x_nF_n\cap E=\varnothing,\qquad
\frac{|(x_nF_n)h\triangle x_nF_n|}{|x_nF_n|}
=\frac{|F_nh\triangle F_n|}{|F_n|}\longrightarrow0,
$$
contradicting \eqref{eq:uniform-bohr-density}. Inversion gives left
syndeticity as well.
\end{proof}

\begin{rem}
The same proof works for every countable amenable group with only
finitely many subgroups of each finite index. The proof of
Lemma~\ref{lem:common-abelian-image} uses only this property.
\end{rem}

For an FC group $H$, each centralizer $C_H(b)=\{h:hb=bh\}$ has
finite index. If $b_1,\ldots,b_r$ generate $H$, then
$Z(H)=\bigcap_j C_H(b_j)$ has finite index. Thus finitely generated
FC groups are virtually abelian; increasing unions and finite
extensions show that every countable virtually FC group is amenable.

\begin{lem}\label{lem:fc-commuting-blocks}
If $H\leq G$ is FC and has finite index, every sequence $(a_n)$ in
$G$ has interval products
$$
b_j=a_{p_j}\cdots a_{q_j}\in H,\qquad q_j<p_{j+1},\qquad
b_ib_j=b_jb_i.
$$
In particular, $\operatorname{FP}_{\uparrow}(b_j)
\subseteq\operatorname{FP}_{\uparrow}(a_n)$.
\end{lem}

\begin{proof}
For a finite-index subgroup $L\leq G$, the partial products
$u_0=e$, $u_r=a_m\cdots a_{m+r-1}$, $1\leq r\leq[G:L]$, give
$[G:L]+1$ elements. Two have the same left coset, and
$$
u_iL=u_jL,\quad i<j
\quad\Longrightarrow\quad
a_{m+i}\cdots a_{m+j-1}=u_i^{-1}u_j\in L.
$$
Start with such a product $b_1\in H$. After choosing
$b_1,\ldots,b_r$, apply this calculation to the tail after $q_r$ and
$$
L_r=\bigcap_{j=1}^rC_H(b_j),\qquad
[G:L_r]\leq[G:H]\prod_{j=1}^r[H:C_H(b_j)]<\infty.
$$
The resulting $b_{r+1}$ commutes with its predecessors. Products of
these disjoint, ordered intervals are increasing products of $(a_n)$.
\end{proof}

\begin{proof}[Proof of Theorem~\ref{thm:bohr}, virtually FC case]
The conclusions pass to supersets, so we may assume
$B=\rho^{-1}(U)$, where $\rho:G\to K$ has compact target and
$U=U^{-1}$ is an open identity neighborhood. Write $b=\|\xi\|^2$, $t=\|P_2\xi\|^2$ and
$c(g)=\Rea\langle\pi(g^2)\xi,\xi\rangle$.
Then $E=\{g\in B:c(g)>t-\varepsilon\}$ is symmetric because
$c(g^{-1})=c(g)$.

Given $(a_n)$, Lemma~\ref{lem:fc-commuting-blocks} gives commuting
interval products $b_j$ with
$\operatorname{FP}_{\uparrow}(b_j)\subseteq
\operatorname{FP}_{\uparrow}(a_n)$. Put $s_j=b_1\cdots b_j$.
Choose an identity neighborhood $V$ with $V^{-1}V\subseteq U$.
A finite cover of $K$ by
left translates of $V$ has a member containing infinitely many
$\rho(s_j)$. Choose $N>1+(b-t)/\varepsilon$ and
$j_1<\cdots<j_N$ in that member. Then
\begin{equation}\label{eq:rec-FC-Bohr-pairs}
\rho(s_{j_r}^{-1}s_{j_s})\in V^{-1}V\subseteq U\qquad(r<s).
\end{equation}
For $v_r=\pi(s_{j_r}^2)\xi$,
$\|v_r\|^2=b$ and $P_2v_r=P_2\xi$. If all
$\Rea\langle v_s,v_r\rangle\leq t-\varepsilon$, then
$$
\begin{aligned}
0\leq\Bigl\|\sum_{r=1}^N(v_r-P_2\xi)\Bigr\|^2
&=N(b-t)+2\sum_{r<s}(\Rea\langle v_s,v_r\rangle-t)\\
&\leq N(b-t)-N(N-1)\varepsilon<0.
\end{aligned}
$$
Hence some $r<s$ gives, using commutation,
$$
\begin{gathered}
g=s_{j_r}^{-1}s_{j_s}=b_{j_r+1}\cdots b_{j_s}
\in\operatorname{FP}_{\uparrow}(a_n),\\
g^2=s_{j_r}^{-2}s_{j_s}^2,\qquad
c(g)=\Rea\langle v_s,v_r\rangle>t-\varepsilon.
\end{gathered}
$$
By \eqref{eq:rec-FC-Bohr-pairs}, $g\in B$, so $E$ is right IP$^*$.
The identity
$(\operatorname{FP}_{\uparrow}(a_n^{-1}))^{-1}
=\operatorname{FP}_{\downarrow}(a_n)$ and $E=E^{-1}$ give left IP$^*$.

If $E$ were not left syndetic, recursively choose
$$
P_n=\{e\}\cup\operatorname{FP}_{\uparrow}(a_1,\ldots,a_n),
\qquad a_{n+1}\notin P_n^{-1}E,
$$
starting with $P_0=\{e\}$. Since $(P_na_{n+1})\cap E=\varnothing$,
every nonempty finite product would avoid $E$, a contradiction.
Thus $K_0E=G$ for a finite nonempty $K_0$, and symmetry gives
$EK_0^{-1}=G$. For every right \Folner{} sequence,
$$
|F_n|\leq\sum_{k\in K_0}|F_nk\cap E|
\leq|K_0|\,|F_n\cap E|+\sum_{k\in K_0}|F_nk\triangle F_n|.
$$
After division by $|F_n|$, the lower density is at least
$1/|K_0|$. Inversion gives the same bound on every left sequence.
\end{proof}

\begin{cor}\label{cor:square-recurrence}
For a measure-preserving action on a probability space $(X,\B,\mu)$
of a countable amenable group $G$ that is finitely generated or
virtually FC, every $A\in\B$,
$\varepsilon>0$ and Bohr identity neighborhood $B$ give a set
$$
\{g\in B:\mu(A\cap T_{g^2}A)>\mu(A)^2-\varepsilon\}
$$
that is left and right syndetic, with one positive lower density
bound for all left and right \Folner{} sequences. It is also left
and right IP$^*$ when $G$ is virtually FC. For every nonnegative,
nonzero $f\in L^2(X,\mu)$,
$$
\lim_n |F_n|^{-1}\sum_{g\in F_n}\int_Xf(f\circ T_{g^2})\,d\mu>0
$$
on every left or right \Folner{} sequence, with a common limit.
In particular, $\delta_A>0$ whenever $\mu(A)>0$.
\end{cor}

\begin{proof}
For the Koopman representation $\pi(g)f=f\circ T_{g^{-1}}$,
constants are fixed by every square. Thus, for $f\geq0$,
$$
t:=\|P_2f\|_2^2\geq\Bigl(\int_Xf\,d\mu\Bigr)^2,\qquad
\langle\pi(g^2)f,f\rangle=\int_Xf(f\circ T_{g^2})\,d\mu\geq0.
$$
For $f=1_A$, the coefficient equals $\mu(A\cap T_{g^2}A)$ and
$t\geq\mu(A)^2$, so Theorem~\ref{thm:bohr} gives the return-set
assertions by inclusion. For $f\ne0$, we have $t>0$.
Theorem~\ref{thm:bohr}, with $B=G$ and $\varepsilon=t/2$, gives
one lower density bound $\eta>0$ for
$E=\{g:\langle\pi(g^2)f,f\rangle>t/2\}$. Therefore
$$
\liminf_n\frac1{|F_n|}\sum_{g\in F_n}\langle\pi(g^2)f,f\rangle
\geq\frac t2\liminf_n\frac{|F_n\cap E|}{|F_n|}
\geq\frac{t\eta}2>0.
$$
Theorem~\ref{thm:square-general} gives the common limit.
\end{proof}

\begin{cor}\label{cor:square-products}
For a countable amenable group $G$ and $D\subseteq G$ with
$d=d_l^*(D)>0$, suppose that $G$ is finitely generated or virtually
FC, or that $d>1/3$. There is $\eta_D>0$ such that every left or
right \Folner{} sequence $F$ admits $B\subseteq G$ with
$$
\overline d_F(B)\geq\eta_D,\qquad BB\subseteq DD^{-1}.
$$
For $d>1/3$, one may take $\eta_D=d(3d-1)/2$.
\end{cor}

\begin{proof}
Apply \cite[Proposition~1.17(1)]{BRM2026} to obtain a
probability-preserving $G$-system and a measurable set $A$ with
$$
\mu(A)=d,\qquad
\{g:\mu(A\cap T_gA)>0\}\subseteq DD^{-1}.
$$
Set $\eta_D=d(3d-1)/2$ if $d>1/3$, and $\eta_D=\delta_A$
otherwise. Corollaries~\ref{cor:large-sets} and
\ref{cor:square-recurrence} give $\delta_A\geq\eta_D>0$.
For any subsequently chosen left or right \Folner{} sequence $F$,
Theorem~\ref{thm:square-general} yields
$$
\lim_n\frac1{|F_n|}\sum_{g\in F_n}\mu(T_g^{-1}A\cap T_gA)
=\delta_A\geq\eta_D.
$$
Now \cite[Theorem~1.24]{BRM2026} gives $B$ with
$$
\overline d_F(B)\geq\eta_D,\qquad
BB\subseteq\{g:\mu(A\cap T_gA)>0\}\subseteq DD^{-1}.
$$
\end{proof}

\section{Limits of square recurrence}\label{sec:recurrence-counterexamples}

A common quaternion construction separates positive square mean,
positive density of near-threshold returns, and the IP$^*$ property.

\begin{lem}\label{lem:quaternion-functions}
For a countably infinite set $I$, write
$$
Q_8=\{\pm1,\pm\mathbf i,\pm\mathbf j,\pm\mathbf k\},\qquad
(\Omega,m)=(Q_8^I,u_{Q_8}^{\otimes I}),
$$
where $u_{Q_8}$ is uniform probability measure. The coordinate functions
$x_i(\omega)=\omega_i$ generate, under pointwise multiplication, a
countable group $G_I=\langle x_i:i\in I\rangle$ that is
locally finite, of nilpotency class two and exponent four. It has a
surjective homomorphism
$$
v:G_I\longrightarrow\bigoplus_{i\in I}\mathbb F_2e_i,\qquad v(x_i)=e_i,
$$
whose kernel $W=Z(G_I)$ has infinite index and consists of
$\{\pm1\}$-valued functions. Every $g\notin W$ has infinitely many
conjugates and is uniformly distributed on $Q_8$. In particular,
$G_I$ is not virtually FC and
\begin{equation}\label{eq:ctr-quaternion-squares}
m\{g^2=1\}=\frac14,\qquad m\{g^2=-1\}=\frac34\qquad(g\notin W).
\end{equation}

On $(X,\mu)=(\Omega\times Q_8,m\otimes u_{Q_8})$, the transformations
and set
$$
T_g(\omega,q)=(\omega,g(\omega)q),\qquad
D=\{1,\mathbf i,\mathbf j,\mathbf k\},\qquad A=\Omega\times D
$$
give a measure-preserving action with $\mu(A)=1/2$ and
\begin{equation}\label{eq:ctr-square-correlation}
\mu(A\cap T_{g^2}A)=\frac12m\{g^2=1\}
=\begin{cases}1/2,&g\in W,\\1/8,&g\notin W.\end{cases}
\end{equation}
For its Koopman representation, $f=1_A-1/2$ belongs to the weakly
mixing subspace and satisfies $P_2f=0$; hence $P_21_A=1/2$ and
$\|P_21_A\|_2^2=1/4$.
\end{lem}

\begin{proof}
For finite $J\subset I$, the group $G_J=\langle x_i:i\in J\rangle$
embeds in $Q_8^{Q_8^J}$, so $G_I=\bigcup_JG_J$ is locally finite.
The identities $a^4=1$ and $[[a,b],c]=1$ hold pointwise.
At $\omega_i=\mathbf i$, $\omega_j=\mathbf j$,
$x_i(\omega)^2=[x_i,x_j](\omega)=-1$; hence the class is two and
the exponent four.

If an identity word has exponent sum $a_i$ at $i$, evaluation
at $\omega_i=\mathbf i$ and all other coordinates $1$ gives
$\mathbf i^{a_i}=1$. Hence exponent sums modulo two define $v$,
with $v(gh)=v(g)+v(h)$ and $v(x_i)=e_i$. Write $\overline a$ for
the image of $a\in Q_8$ in $Q_8/\{\pm1\}\cong\mathbb F_2^2$. Then
\begin{equation}\label{eq:ctr-quaternion-quotient}
\overline{g(\omega)}=\sum_i v(g)_i\overline{\omega_i}.
\end{equation}
Thus $W=\ker v$ consists of $\{\pm1\}$-valued functions and
$W\subseteq Z(G_I)$, while $[G_I:W]=\infty$.
If $v(g)_i=1$ and $j,k$ are distinct coordinates unused in a word
for $g$, evaluate at $\omega_i=\mathbf i$, $\omega_j=\mathbf j$,
and all other coordinates $1$. This gives
$$
g(\omega)=\pm\mathbf i,\qquad
[g,x_j](\omega)=-1,\qquad[g,x_k](\omega)=1.
$$
The commutators with the unused coordinates are pairwise distinct.
Consequently $W=Z(G_I)$ and every $g\notin W$ has infinitely many
conjugates. If an FC subgroup $L$ had finite index, a left-coset
decomposition would give, for $h\in L$,
$$
|\{ghg^{-1}:g\in G_I\}|
\leq[G_I:L]\,|\{\ell h\ell^{-1}:\ell\in L\}|<\infty.
$$
Then $L\subseteq W$ and $[G_I:W]\leq[G_I:L]<\infty$, a contradiction.

For $g\notin W$, the right side of
\eqref{eq:ctr-quaternion-quotient} is uniform on $\mathbb F_2^2$.
Replacing a coordinate with $v(g)_i=1$ by
$-\omega_i$ preserves $m$ and sends $g$ to $-g$.
Therefore
$$
m\{g=a\}=m\{g=-a\}
=\tfrac12m\{g\in\{a,-a\}\}=\tfrac18\qquad(a\in Q_8).
$$
Only $\pm1$ square to $1$, proving \eqref{eq:ctr-quaternion-squares}.

Each $g$ uses finitely many coordinates, so $T_g$ is measurable,
and $T_gT_h=T_{gh}$. Since $q\mapsto g(\omega)q$ permutes $Q_8$,
each $T_g$ preserves $m\otimes u_{Q_8}$.
The partition $D\sqcup(-D)=Q_8$ gives $\mu(A)=1/2$ and
$u_{Q_8}(D\cap g(\omega)^2D)=\frac12 1_{\{g^2=1\}}(\omega)$.
Integration proves \eqref{eq:ctr-square-correlation}.

The closed subspace
$$
\HH_-=\{F\in L^2(X,\mu):F(\omega,-q)=-F(\omega,q)\text{ a.e.}\}
$$
contains $f$ and is invariant under $\pi(g)$ and $\pi(g)^*$.
On this subspace,
$$
\pi(x_i^2)F=\epsilon_iF,\qquad
\epsilon_i(\omega)=\omega_i^2\in\{\pm1\},\qquad
m\{\epsilon_i=1\}=\tfrac14.
$$
Independence gives, for an enumeration $I=\{i_1,i_2,\ldots\}$
and any $(s_i)\in\{\pm1\}^I$,
\begin{equation}\label{eq:ctr-no-atoms}
m\{\epsilon_i=s_i\ (i\in I)\}
\leq\prod_{j=1}^r m\{\epsilon_{i_j}=s_{i_j}\}
\leq(3/4)^r\longrightarrow0.
\end{equation}
A nonzero finite-dimensional invariant subspace of $\HH_-$ would
contain a common eigenvector $F\ne0$ of the commuting involutions
$\pi(x_i^2)$. The equations $\epsilon_iF=s_iF$ hold on one conull
set, so $F$ is supported on
$\bigcap_i\{\epsilon_i=s_i\}\times Q_8$, a null set by
\eqref{eq:ctr-no-atoms}. This is a contradiction. The projection onto $\HH_-$
commutes with $\pi(G_I)$, so its image of any finite-dimensional
invariant subspace is invariant and must be zero. Thus $\HH_-$
lies in the weakly mixing subspace. Taking all $s_i=1$ also excludes
nonzero vectors of $\HH_-$ fixed by all squares, so $P_2\HH_-=0$.
\end{proof}

\begin{proof}[Proof of Theorem~\ref{thm:no-fc-counterexample}]
Take Lemma~\ref{lem:quaternion-functions} with $I=\mathbb N$.
By \eqref{eq:ctr-square-correlation}, the return set is $W$.
Since $[G:W]=\infty$, the subgroup $W$ is neither left nor right
syndetic. Lemma~\ref{lem:subgroup-density} gives
$|F_n\cap W|/|F_n|\to0$ on every left or right \Folner{} sequence.
Both product orders avoid $W$, because
$$
v(x_{i_1}\cdots x_{i_r})=v(x_{i_r}\cdots x_{i_1})
=e_{i_1}+\cdots+e_{i_r}\ne0\qquad(i_1<\cdots<i_r).
$$
Thus $W$ is neither right nor left IP$^*$. Lemma~\ref{lem:quaternion-functions}
places $f=1_A-1/2$ in the weakly mixing subspace, and
\eqref{eq:ctr-square-correlation} gives
$$
\begin{gathered}
\langle\pi(g^2)f,f\rangle
=\mu(A\cap T_{g^2}A)-\tfrac14
=\begin{cases}1/4,&g\in W,\\-1/8,&g\notin W,\end{cases}\\
\{g:\langle\pi(g^2)f,f\rangle>-\delta\}=W
\qquad(0<\delta<1/8).
\end{gathered}
$$
\end{proof}

For $G_n=\langle x_1,\ldots,x_n\rangle$, we have
$[G_n:G_n\cap W]=2^n$ and $G=\bigcup_nG_n$. Thus $(G_n)$ is two-sided
\Folner{}, and the positive square mean is explicit:
$$
|G_n|^{-1}\sum_{g\in G_n}\mu(A\cap T_{g^2}A)
=\frac18+\frac38\,2^{-n}.
$$

\begin{thm}\label{thm:two-generated-square-ip-counterexample}
There is a two-generated solvable group $H$ of derived length at
most three with a weakly mixing measure-preserving action on a
probability space $(X,\B,\mu)$, a set $A$ of measure $1/2$ and a
sequence $(a_n)$ such that
$$
\mu(A\cap T_{g^2}A)=\frac18
\qquad
\bigl(g\in\operatorname{FP}_{\uparrow}(a_n)
\cup\operatorname{FP}_{\downarrow}(a_n)\bigr).
$$
For every $0<\varepsilon<1/8$, the return set
$$
E_\varepsilon=\{h\in H:\mu(A\cap T_{h^2}A)>\mu(A)^2-\varepsilon\}
$$
is neither left nor right IP$^*$ but has density one along every
left and every right \Folner{} sequence.
\end{thm}

\begin{proof}
Use Lemma~\ref{lem:quaternion-functions} with $I=\mathbb Z$ and
$N=G_I$. The shift gives
$$
(R\omega)_i=\omega_{i-1},\qquad
\sigma(n)=n\circ R^{-1},\qquad
H=N\rtimes_\sigma\langle t\rangle,
\qquad tnt^{-1}=\sigma(n).
$$
Since $\sigma(x_i)=x_{i+1}$, we have $H=\langle x_0,t\rangle$,
$H'\subseteq N$ and $H'''\subseteq N''=\{e\}$.
The group is amenable as an extension of the locally finite group
$N$ by $\mathbb Z$. Its measure-preserving action extends by
$T_t(\omega,q)=(R\omega,q)$, since
$$
T_tT_nT_t^{-1}(\omega,q)
=(\omega,n(R^{-1}\omega)q)=T_{\sigma(n)}(\omega,q).
$$
Retain $A$ and $f$ from Lemma~\ref{lem:quaternion-functions} and
take $a_i=x_i$, $i\geq1$. The two finite-product sets lie in
$N\setminus W$, so \eqref{eq:ctr-square-correlation} gives the
value $1/8$ and both IP$^*$ failures. For the projections onto vectors fixed by all squares,
$P_2^H\leq P_2^N$ and Lemma~\ref{lem:quaternion-functions}
gives $P_2^Nf=0$; hence $P_2^Hf=0$.
The centered correlation is $-1/8$ on both finite-product sets.

To prove weak mixing, take a function $F\in L^2(X\times X)$
invariant under the diagonal action. For each $(a,b)\in Q_8^2$,
$$
(m\otimes m)\{(\omega_i,\omega_i')\ne(a,b)\text{ for }1\leq i\leq r\}
=(1-1/64)^r\longrightarrow0.
$$
Thus almost every $(\omega,\omega')$ contains every pair $(a,b)$
among its coordinates, so $n\mapsto(n(\omega),n(\omega'))$
maps $N$ onto $Q_8^2$.
Since $N$ is countable and the fibers are finite, Fubini puts all
invariance equations on one conull set of base pairs. On this set,
invariance and surjectivity give, for $n\in N$ and $q,q'\in Q_8$,
$$
\begin{aligned}
F(\omega,n(\omega)q;\omega',n(\omega')q')
&=F(\omega,q;\omega',q'),\\
F(\omega,q;\omega',q')&=F(\omega,1;\omega',1)
=:\varphi(\omega,\omega').
\end{aligned}
$$
The remaining invariance says $\varphi\circ(R\times R)=\varphi$.
This is the Bernoulli shift on $(Q_8^2)^{\mathbb Z}$, which is
ergodic, so $F$ is constant.

Finally, $W=\ker v$ is $\sigma$-invariant, and
$$
(nt^k)^2=n\sigma^k(n)t^{2k},\qquad
v(n\sigma^k(n))=\sum_i v(n)_i(e_i+e_{i+k}).
$$
For $k\ne0$ and $n\notin W$, the sum is nonzero, since a nonzero
finitely supported vector cannot be periodic under a nonzero shift.
Lemma~\ref{lem:quaternion-functions} makes $n\sigma^k(n)$ uniform
on $Q_8$. For $h\notin N\cup(W\rtimes\langle t\rangle)$,
$T_tA=A$ therefore gives
$$
\mu(A\cap T_{h^2}A)
=\frac18\sum_{a\in Q_8}u_{Q_8}(D\cap aD)
=u_{Q_8}(D)^2=\frac14.
$$
Both subgroups have infinite index:
$[H:N]=\infty$ and
$[H:W\rtimes\langle t\rangle]=[N:W]=\infty$.
Lemma~\ref{lem:subgroup-density} now gives, for every left or
right \Folner{} sequence,
\begin{equation*}
0\leq1-\frac{|F_n\cap E_\varepsilon|}{|F_n|}
\leq\frac{|F_n\cap N|+|F_n\cap(W\rtimes\langle t\rangle)|}{|F_n|}
\longrightarrow0.\qedhere
\end{equation*}
\end{proof}

\section{Examples}\label{sec:examples}

\begin{exam}\label{ex:flip}
The antidiagonal limit need not be positive or a projection,
even for a finite group.
\end{exam}

\begin{proof}
For the standard irreducible representation
$\pi:S_3\to\U(\C^2)$, take
$U_g=\pi(g)\otimes I$, $V_g=I\otimes\pi(g)$ and $F_n=S_3$.
Equation~\eqref{eq:op-Schur} gives
$$
D_{U,V}=\frac1{|S_3|}\sum_{g\in S_3}
\pi(g)^*\otimes\pi(g)=\frac12\mathsf F_2,
\qquad \spec(D_{U,V})=\{-1/2,1/2\}.
$$
\end{proof}

\begin{exam}\label{ex:strong}\label{ex:square-strong}
For a countable locally finite group, antidiagonal and square
averages can fail to converge strongly along the same left
\Folner{} sequence. This occurs for unitary representations and
measure-preserving actions. The adjoints of square averages can
fail to converge strongly along a right \Folner{} sequence.
\end{exam}

\begin{proof}
On $\N_0=\{0,1,\ldots\}$, let $K_n$ be the finite permutation
group fixing every $j>n$. Define
$$
G=\bigcup_{n\geq1}K_n,\qquad
t_n=(0,n+1,n+2),\qquad F_n=K_nt_n.
$$
The increasing union makes $G$ countable and locally finite.
Every $h\in G$ eventually belongs to $K_n$, and then
$hF_n=F_n$; hence $(F_n)$ is left \Folner{}.
For $g=kt_n\in F_n$,
$$
g(0)=n+1,\qquad g(n+1)=n+2,\qquad g^2(0)=n+2.
$$
Use $\pi(g)e_j=e_{g(j)}$ on $\ell^2(\N_0)$, and take
$(U,V)=(I,\pi)$ for the antidiagonal averages. We obtain
$$
D_{F_n}e_0=e_{n+1},\qquad
A_{F_n}e_0=e_{n+2},\qquad
A_{F_n^{-1}}^*e_0=e_{n+2}.
$$
Since $\|e_j-e_k\|_2=\sqrt2$ for $j\ne k$, none of these
sequences converges in norm; $(F_n^{-1})$ is right \Folner{}.
Theorem~\ref{thm:main} also gives
$$
\ran D_{I,\pi}\subseteq
\{\xi\in\ell^2(\N_0):\pi(g)\xi=\xi\text{ for every }g\in G\}
=\{0\},
$$
because an invariant vector has constant coordinates.

For measure-preserving actions, equip $X=\{-1,1\}^{\N_0}$ with
fair product measure and set $(S_gx)_j=x_{g^{-1}(j)}$.
Its Koopman representation satisfies, for $r_j(x)=x_j$,
$$
V_gr_j=r_{g(j)},\qquad
\Bigl\|\sum_j a_jr_j\Bigr\|_2^2=\sum_j|a_j|^2
$$
for every finitely supported family $(a_j)$. Taking the trivial
action for $U$, the three identities above hold with $r_j$ in
place of $e_j$ and give the same failure of norm convergence.
\end{proof}

General results on one-sided mean ergodic averages appear in
\cite{BF}.

\begin{exam}\label{ex:antidiagonal-below-square}
Two commuting, individually ergodic actions of a finite group
can satisfy
$$
\frac1{|G|}\sum_{g\in G}\mu(T_gA\cap S_gA)
=\frac5{54}<\frac19=\mu(A)^2.
$$
\end{exam}

\begin{proof}
Equip $X=G=S_3\times C_3$ with uniform probability measure.
The actions $T_gx=gx$, $S_gx=xg^{-1}$ commute and are transitive,
hence individually ergodic. Writing $C_3$ additively, take
$$
A=\bigl((\{e\}\cup\mathcal T)\times\{0\}\bigr)
\cup\{(r,0),(r,1)\},\qquad
r=(123),\qquad \mu(A)=\frac13,
$$
where $\mathcal T$ is the set of three transpositions.

For $a=(s,u),b=(t,v)\in A$, the equivalence
$ga=bg^{-1}\iff(ga)^2=ba$ and bijectivity of squaring on $C_3$ give
$$
|\{g\in G:ga=bg^{-1}\}|=|\{h\in S_3:h^2=ts\}|
=\begin{cases}
4,&ts=e,\\
1,&ts\in\{r,r^2\},\\
0,&ts\in\mathcal T.
\end{cases}
$$
Among the even permutations, $e,r,r^2$ occur in $A$ with
multiplicities $1,2,0$; each transposition occurs once.
Pairs of opposite parity contribute zero. Therefore
\begin{equation}\label{eq:sq-square-root-count}
\begin{aligned}
\sum_{g\in G}|gA\cap Ag^{-1}|
&=\sum_{a,b\in A}|\{g\in G:ga=bg^{-1}\}|\\
&=\underbrace{4\cdot1^2+2\cdot1\cdot2+2^2}_{s,t\text{ even}}
+\underbrace{3\cdot4+6\cdot1}_{s,t\text{ odd}}
=30.
\end{aligned}
\end{equation}
By \eqref{eq:sq-square-root-count},
$$
\frac1{|G|}\sum_{g\in G}\mu(T_gA\cap S_gA)
=\frac{30}{|G|^2}=\frac{30}{18^2}=\frac5{54}.
$$
\end{proof}

The actions extend to $G\times\Z$ by a trivial $\Z$-action.
They remain individually ergodic, and averaging over
$G\times[-n,n]$ still gives $5/54$.

\end{document}